\documentclass[12pt]{article}

\usepackage{amssymb}
\usepackage{amsmath}
\usepackage{latexsym}
\newtheorem{theorem}{Theorem}

\newtheorem{lemma}[theorem]{Lemma}

\newtheorem{definition}[theorem]{Definition}
\newtheorem{remark}[theorem]{Remark}

\newcommand{\nl}{\newline}

\newcommand{\dom}{{\rm Dom}}

\newcommand{\kernel}{{\rm Ker}}

\newcommand{\N}{{\bf N}}

\newcommand{\R}{{\bf R}}

\newcommand{\cC}{{\cal C}}

\newcommand{\cW}{{\cal W}}

\newcommand{\parder}[2]{\frac{\partial{#1}}{\partial{#2}}}

\newcommand{\inprod}[2]{{\langle{#1},{#2}\rangle}}

\newcommand{\half}{\frac{1}{2}}

 \newcommand{\ia}{({\rm i})}
 \newcommand{\ib}{({\rm ii})}
 \newcommand{\ic}{({\rm iii})}

\numberwithin{equation}{section}

\title{Spectral stability estimates for elliptic operators subject to domain transformations with non-uniformly bounded
gradients }
\author{Gerassimos Barbatis\footnote{Corresponding author}\ \ and Pier Domenico Lamberti}
\date{}
\begin{document}
\maketitle

\begin{abstract}
\noindent
We consider  uniformly elliptic operators with Dirichlet or Neumann homogeneous boundary conditions on a domain $\Omega $ in ${\mathbb{R}}^N$. We consider deformations $\phi (\Omega )$ of $\Omega $ obtained by means of a locally Lipschitz homeomorphism $\phi $ and we estimate the variation  of the eigenfunctions and eigenvalues upon variation of $\phi $. We prove general stability estimates without using  uniform upper bounds  for the gradients of the maps $\phi$. As an application, we obtain estimates on the rate of convergence for eigenvalues and eigenfunctions when a domain with an outward cusp is approximated by a sequence of Lipschitz domains.
\end{abstract}

\vspace{11pt}

\noindent
{\bf Keywords:} elliptic equations, spectral stability, domain perturbation, cusps.

\vspace{6pt}
\noindent
{\bf 2010 Mathematics Subject Classification:} 35J25, 47A75, 47B25.

\section{Introduction}

Let $\Omega $ be a bounded domain ({\it i.e.}, a bounded connected open set) in $\R^N$ and $\phi $ a locally Lipschitz homeomorphism between  $\Omega $ and another
bounded domain $\phi (\Omega )$ in $\R^N$.
For fixed real coefficients $A_{ij}$ defined in the whole of $\R^N$ with $A_{ij}=A_{ji}$ and  satisfying the uniform ellipticity condition (\ref{ellip}),  we consider in $\phi (\Omega )$ the operator $L$ defined formally by
\begin{equation}
\label{classic}
-\sum_{i,j=1}^N\frac{\partial}{\partial x_i }\left( A_{ij}(x)\frac{\partial u}{\partial x_j} \right) ,
\end{equation}
and subject  to the Dirichlet boundary condition $u=0$ on  $\partial \phi (\Omega )$  or the Neumann boundary condition
$$\sum_{i,j=1}^N A_{ij}\frac{\partial u}{\partial x_j}\nu_{i}=0,\ \  {\rm on}\ \partial \phi (\Omega ),$$
where $\nu=(\nu_1,\dots ,\nu_N)$ denotes the outer unit normal to $\partial\phi(\Omega)$. Following the approach developed in \cite{babula},
we prove  estimates for the deviation of the eigenvalues and eigenfunctions corresponding to the domain $\phi (\Omega )$ from those corresponding to a perturbation $\tilde \phi (\Omega )$ of $\phi (\Omega )$. In particular, we improve the results of \cite{babula} in two respects: we provide estimates
which allow dealing with  possibly singular maps $\phi $  and we improve the exponents appearing in the appropriate measures of vicinity of $\phi (\Omega )$ and $\tilde \phi (\Omega )$.

The regular case of  globally Lipschitz homeomorphisms  $\phi$ was investigated in \cite{babula} where  estimates for the variation of the resolvents, eigenfunctions and eigenvalues were proved under the assumption that the gradients of the maps $\phi$ and their inverses have a uniform  upper bound. Those estimates can be applied for example to the case of uniform families of domains with Lipschitz continuous boundaries.
However, if $\phi (\Omega )$ has  boundary degenerations stronger than those of $\Omega$ (for example, $\Omega$ has a Lipschitz continuous boundary  while $\phi (\Omega )$ has a cusp at the boundary)  one cannot assume that $\phi $ has a bounded gradient. This problem might  be overcome
by approximating  $\phi (\Omega )$ by means of suitable domains $\phi_{\epsilon }(\Omega )$, $\epsilon >0$, where $\phi_{\epsilon }$ are globally Lipschitz continuous maps: one would find estimates depending on $\epsilon $ and eventually would pass to the limit as $\epsilon \to 0$.  However, the gradients of the maps $\phi_{\epsilon }$  would not necessarily have a uniform upper bound, hence the results of \cite{babula} could not be used in  this limiting procedure. Thus, it is desirable to prove  stability estimates independent of  $\|\nabla \phi\|_{L^\infty(\Omega )}$.
In this paper, we prove general stability estimates without using any uniform  upper bound for $\|\nabla \phi\|_{L^\infty(\Omega )}$.
These estimates are expressed in terms of a certain measure of vicinity  $\delta _{q}(\phi ,\tilde \phi )$   of $\phi $ and $\tilde \phi$ which  reduces
 to the Sobolev norm $\| \phi -\tilde\phi \|_{W^{1,q}(\Omega )}$ in regular cases,  see (\ref{vicinity}) for the precise definition and Remark~\ref{sobo}.

 Similarly to \cite{babula} the estimates for the variation of eigenvalues and eigenfunctions are deduced from corresponding estimates for
the variation of   resolvent operators in the Hilbert-Schmidt class ${\mathcal{ C}}^2 $. Note that the resolvent $(L+1)^{-1}$ of the operator $L$ belongs to the Schatten class ${\mathcal{C}}^r$, $1\le r <\infty $  if and only if the eigenvalues $\lambda_n$ of $L$ satisfy
 $$
 \sum_{n=1}^{\infty }\frac{1}{(\lambda_n+1)^r} <\infty ,
 $$
and,  in the case of smooth domains, this holds provided $r>N/2$. Condition $r>N/2$ is used in \cite{babula} and turns out to spoil the exponents in the  stability estimates. If one is interested only in eigenvalues and eigenfunctions (and not in the solutions to the Poisson problem $Lu=f$), it is convenient to replace the resolvent $(L+1)^{-1}$ by suitable powers of it. Indeed,
the operator $(L+1)^{-k}$ belongs to any fixed  Schatten class ${\mathcal{C}}^r$ provided $k\in \N$ is large enough. The power $k$ plays  no essential  role in the estimates for eigenvalues and eigenfunctions:  this simple but crucial observation enables us to improve the estimates of \cite{babula}.\\
In the case of transformations $\phi$ with uniformly bounded gradients considered in \cite{babula}, the new estimate for eigenvalues reads

\begin{equation}
\left( \sum_{n=1}^{\infty}\left[
(\tilde \lambda_n +1)^{-k} -  (\lambda_n +1)^{-k}\right]^2\right)^{1/2}
\leq c \delta_{\frac{2q_0}{q_0-2}}(\phi ,\tilde\phi) ,
\label{res2'intro}
\end{equation}
where $ \lambda_n$, $\tilde\lambda_n$ denote the eigenvalues in $\phi (\Omega )$ and $\tilde \phi (\Omega )$ respectively,  see Theorem~\ref{thm:series}.  Here $q_0\in ]2,\infty ]$
is a suitable parameter related to a summability assumption (property (P)) on the eigenfunctions and their gradients, see Definition~\ref{growth}. It turns out that in the case of sufficiently smooth domains (say, of class $C^{1,1}$) and sufficiently smooth coefficients $A_{ij}$ (say, Lipschitz continuous), one can take $q_0=\infty$, hence the estimate (\ref{res2'intro}) is expressed in terms of $\delta _2(\phi ,\tilde \phi )$. (Note that in \cite{babula} the best measure of vicinity appearing in the estimates is $\delta_{N+\epsilon }(\phi ,\tilde\phi) $, for any $\epsilon >0$ and it is much worse than $\delta_2(\phi, \tilde\phi )$ if $N>2$). In the general case of possibly singular transformations $\phi$, the term $\delta_{2q_0/(q_0-2)}(\phi ,\tilde\phi)$ in (\ref{res2'intro}) has to be replaced by $(1+\delta_s(\phi ,\tilde\phi ))\delta_{2q_0/(q_0-2)}(\phi ,\tilde\phi)$ for a suitable $s\geq 1$, where the extra summand  appears only for technical reasons and is not important for applications.

Estimate (\ref{res2'intro}) is first applied to uniform families of domains with Lipschitz continuous boundaries as in \cite{babula}. In this case, the construction of
appropriate transformations $\phi$ leads to the estimate
\begin{equation}\label{measintro}
\left( \sum_{n=1}^{\infty}\left[
(\lambda_n[\Omega_1] +1)^{-k} -  (\lambda_n[\Omega_2] +1)^{-k}\right]^2\right)^{1/2}
\leq c |\Omega_1\vartriangle \Omega_2|^{\frac{1}{2}-\frac{1}{q_0}},
\end{equation}
provided $\Omega_1 $ and $\Omega_2$ belong to the same Lipschitz class and the Lebesgue measure $|\Omega_1\vartriangle \Omega_2|$ of the symmetric difference of $\Omega_1$ and $\Omega_2$  is small enough, see Theorem~\ref{thm:C11}.
Analogous estimates for the variation of the eigenfunctions are proved in Theorems~\ref{thm:eigen}, \ref{thm:C11}.

We then apply our general stability estimates to the case of the Dirichlet Laplacian on  a domain $\Omega$ with an exterior power-type cusp of exponent $\alpha $ sufficiently close to 1 (the case $\alpha =1$ is clearly the regular Lipschitz case). We approximate $\Omega$ by a sequence $\Omega_{\delta }$, $\delta >0$, of domains with Lipschitz continuous boundaries and estimate the rate of convergence of the eigenvalues and eigenfunctions in terms of  $ |\Omega \setminus \Omega_{\delta }  |^{b(\alpha )}  $, where $0<b(\alpha )<1$ is an explicit exponent depending only on $N$ and $\alpha$ (with $b(1)=1/2$ as expected from (\ref{measintro})). To do so, we establish the validity of property (P) in the domains $\Omega_{\delta }$ by means of an apriori estimate of Maz'ya and Plamenevskii~\cite[p.~4]{mp} and a bootstrap argument, see Theorems~\ref{bootstrap_thm} and \ref{pcusp}. According to the strategy explained above, we then construct suitable maps
 $\phi_{\epsilon }:\Omega_{\delta}\to \Omega_{\epsilon}$  and get  estimates in terms of $|\Omega _{\epsilon }\setminus \Omega_{\delta}|^{b(\alpha )}$. By letting $\epsilon \to 0$ we obtain the desired estimate.

We note that  in the case of suitable uniform families of domains with Lipschitz continuous boundaries it was proved in \cite{buladir, bulaneu} that
\begin{equation}
\label{measintrobur}
|\lambda_n[\Omega_1]-\lambda_n[\Omega_2]|\le c_n |\Omega_1\triangle\Omega_2|^{1-\frac{2}{q_0}},
\end{equation}
where $q_0 $ is as above. Moreover, in \cite{highsharp} it is proved  that the exponent $1-2/q_0$ is sharp, see also \cite{lape}.
Clearly, in estimate (\ref{measintro}) we do not obtain the sharp exponent.
The fact that our exponent is exactly twice the sharp one seems to indicate that a
variation of our method could lead to the optimal exponent.  However, we note that our method has the advantage of providing stability estimates also for eigenfunctions and large enough powers of the resolvents. Such estimates cannot be obtained by the methods of \cite{buda} and \cite{buladir, bulaneu} which make use of the variational characterization of the eigenvalues. We note that while stability estimates for eigenvalues have been extensively studied in recent years, the corresponding problem for eigefunctions is much less investigated. In this  respect  we mention the article of Pang \cite{pang} where probabilistic methods are used to obtain a stability estimate for the ground state of the Dirichlet Laplacian on a simply connected planar domain.

We note that estimates of the type (\ref{measintrobur}) have been recently obtained by Lemenant and Milakis~\cite{mile} for the first eigenvalue of the Dirichlet Laplacian in Reifenberg flat domains. We also note that stability estimates for the eigenvalues of uniformly elliptic operators with Dirichlet or Neumann boundary conditions on domains with  continuous boundaries were proved in Burenkov and Davies~\cite{buda} and in \cite{highsharp, bulahigh} where the vicinity of the domains is expressed in terms of a  variant of the Hausdorff distance.
For more references on this subject we refer to \cite{babula} and to the survey paper \cite{bulalanz}.

This paper is organized as follows. In Section~\ref{operators} we set the problem. In Section~\ref{perturb} we prove  stability estimates for resolvents, eigenvalues and eigenfunctions in terms of $\delta_q(\phi,\tilde\phi )$. In Section~\ref{applicationsec} we discuss some applications to domains with Lipschitz continuous boundaries as well as to domains with power-type cusps at the boundary, and we prove estimates in terms of the Lebesgue measure.

\section{Elliptic operators and singular domain transformations}
\label{operators}

Let $\Omega$ be an arbitrary  bounded domain in $\R^N$. We
consider a family of
domains $\phi\left(\Omega\right)$ in $\R^N$
parametrized by locally Lipschitz homeomorphisms $\phi$ of $\Omega $ onto $\phi (\Omega )$.
More precisely, we consider the  family of transformations \begin{eqnarray}
\label{phis}
& &\Phi (\Omega) := \left\{
\phi\in\left(W^{1,\infty}_{loc}(\Omega)\cap L^{\infty }(\Omega )   \right)^{N}:\,
{\mathrm{the\ continuous\ representative\ of}}\ \phi \right.\nonumber
\\
& &\left.
\qquad\qquad
\qquad\qquad\qquad\qquad
{\mbox{is injective and}}\
\phi^{(-1)}\in\left(W^{1,\infty}_{loc}(\phi (\Omega ))\right)^{N}
\right\},
\end{eqnarray}
where $W^{1,\infty}_{loc}(\Omega)$ denotes the Sobolev space
of the functions in $L^{\infty }_{loc}\left(\Omega\right)$
which have weak derivatives of first order in
$L^{\infty}_{loc}\left(\Omega\right)$. Observe that if $\phi\in \Phi (\Omega )$ then $\phi$ is locally Lipschitz continuous.
Note also that if $\phi \in \Phi (\Omega )$ then
$\phi (\Omega)$ is also a bounded domain.
Moreover, any transformation $\phi \in \Phi (\Omega )$ allows changing variables in integrals in the standard way.

Let $A=(A_{ij})_{i,j=1, \dots , N}$ be a real symmetric matrix-valued function defined on $\R^N$ such that $A_{ij}\in L^{\infty }(\R^N )$ for all $i,j=1, \dots , N$ and
\begin{equation}
\label{ellip}
\theta^{-1} |\xi |^2 \le \sum_{i,j=1}^NA_{ij}(x)\xi_i\xi_j \le \theta |\xi |^2 ,
\end{equation}
for all $x, \xi \in \R^N$ and some $\theta\geq 1$. This matrix will be fixed throughout the paper.

Let $\phi\in \Phi(\Omega)$ and let $\cW$ denote either $W^{1,2}_0(\phi(\Omega))$ or $W^{1,2}(\phi(\Omega))$. Here $W^{1,2}(\phi(\Omega))$ denotes the standard Sobolev space of functions in $L^2(\phi (\Omega ))$ with first order weak derivatives in $L^2(\phi (\Omega ))$ endowed with its usual norm, and
  $W^{1,2}_0(\phi(\Omega))$ denotes the closure in $W^{1,2}(\phi(\Omega))$ of the $C^{\infty }$-functions with compact support in $\Omega$. We consider a non-negative self-adjoint operator $L$ on $L^2(\phi(\Omega))$ given formally by (\ref{classic}) and satisfying Dirichlet or Neumann boundary conditions on $\partial\phi (\Omega )$. More precisely, $L$ is defined as the self-adjoint operator on $L^2(\phi(\Omega))$ canonically associated with the quadratic form $Q_L$ given by
\begin{equation}
\label{quadratic}
\dom(Q_L)=\cW \, , \quad Q_L(v)=\int_{\phi(\Omega)} \sum_{i,j=1}^NA_{ij}(y)\parder{v}{y_i}\parder{\bar v}{y_j}
   dy ,
\end{equation}
for all $v\in \cW$.
We now consider the operator $H$ on $L^2(\Omega)$ obtained by pulling-back $L$ to $\Omega$ as follows.
Let $C_{\phi }$ be the operator from $L^2(\phi (\Omega ))$ to $L^2(\Omega )$ defined by $C_{\phi }v=v\circ \phi $ for all
$v\in L^2(\phi (\Omega ))$.
Let $v\in W^{1,2}(\phi (\Omega ))$ be given and let $u= C_{\phi }v$. Observe that
\[
\int_{\phi(\Omega)}|v|^2dy =\int_{\Omega}|u|^2|\det\nabla\phi(x)| \, dx\; .
\]
Moreover a simple computation shows that
\[
\int_{\phi(\Omega)}   \sum_{i,j=1}^NA_{ij}(y)\parder{v}{y_i}\parder{\bar v}{y_j} dy =\int_{\Omega }\sum_{i,j=1}^Na_{ij}(x)\parder{u}{x_i}\parder{\bar u}{x_j}
|\det\nabla\phi(x)| \, dx\; ,
\]
where $a=(a_{ij})_{i,j=1, \dots , N}$ is the matrix valued function defined on $\Omega$ by
\begin{eqnarray*}
a_{ij}&=&\sum_{r,s=1}^N \Big(A_{rs}  \parder{  \phi_i^{(-1)}}{y_r}   \parder{\phi_j^{(-1)}}{y_s}\Big)\circ \phi\\
&=& ((\nabla \phi)^{-1} A(\phi)  (\nabla\phi)^{-t})_{ij} \; .
\end{eqnarray*}
Here $(\nabla \phi )^{-t}$ denotes the transpose of the inverse of the matrix $\nabla \phi$.
The operator $H$ is defined as the non-negative self-adjoint operator on the Hilbert space $L^2(\Omega , |\det\nabla\phi(x)|\, dx)$
associated with the closure of the quadratic form $Q_H$ with $\dom(Q_H  )=C_{\phi}[\cW]$ and
\[
Q_{H}(u)=\int_{\Omega}\sum_{i,j=1}^Na_{ij}\parder{u}{x_i}\parder{\bar {u}}{x_j}|\det\nabla\phi(x)|\, dx  , \qquad u\in \dom(Q_H).
\]
We note that $H$ is not necessarily uniformly elliptic. We also note that, equivalently, $H$ can be defined as
$$
H=C_{\phi }L C_{\phi ^{(-1)}}.
$$
In particular $H$ and $L$ are unitarily equivalent and the operator $H$ has compact resolvent if and only if $L$ has compact resolvent.
We set $$g(x):=|\det\nabla\phi(x)|,$$ for all $x\in\Omega$, and we denote by $\inprod{\cdot}{\cdot}_g$
the inner product in $L^2(\Omega,g\, dx)$ and also in $(L^2(\Omega,g\, dx))^N$.

\
\section{Stability estimates}
\label{perturb}

In this section we shall consider maps $\phi$ with the properties described in Section \ref{operators}, and we make the additional assumption that $\phi$ and its inverse $\phi^{(-1)}$ are Lipschitz continuous. We note that in this case
\[
C_{\phi}[W^{1,2}(\phi(\Omega))] =W^{1,2}(\Omega)   \;\; \mbox{ and }\;\;  C_{\phi}[W^{1,2}_0(\phi(\Omega))] =W^{1,2}_0(\Omega) .
\]
In this context we give an additional definition. We define $T:L^2(\Omega, \, g\, dx)\to (L^2(\Omega,\, g\, dx))^N$ to be the operator with domain $\dom(T)=C_{\phi}[\cW]$ and $Tu=a^{1/2}\nabla u$. We then have
\begin{equation}
H=T^{(*)_g}T .
\label{t}
\end{equation}
Here the adjoint $T^{(*)_g}$ of $T$ is understood with respect to the inner product of $L^2(\Omega, \, g\, dx)$ and this has been
emphasized in the notation. However, in the sequel we shall simply write $T^{*}$ instead of  $T^{(*)_g}$, unless it is necessary to distinguish  two different scalar products.

Let $\phi$ and $\tilde\phi$ be two such maps on $\Omega$ and let $L$ and $\tilde L$ be the corresponding operators on $\phi(\Omega)$ and $\tilde\phi(\Omega)$ defined as in Section
\ref{operators}. We assume that either $L$ and $\tilde L$ both satisfy Dirichlet boundary conditions or $L$ and $\tilde L$ both satisfy Neumann boundary conditions.
We shall use a tilde to distinguish the objects corresponding to $L$ from those corresponding to $\tilde L$.
Our aim is to compare $L$ and $\tilde L$ and to do this we shall compare the respective pull-backs $H$ and $\tilde H$.
Since $H$ and $\tilde H$ act on different Hilbert spaces -- $L^2(\Omega ,  g\, dx)$ and $L^2(\Omega , \tilde g\, dx)$ -- we shall use the canonical unitary operator,
\[
w: L^2(\Omega , g\, dx) \longrightarrow L^2(\Omega , \tilde g\, dx) \;\; , \qquad u\mapsto wu\, ,
\]
defined as the multiplication by the function  $w:=g^{1/2}{\tilde g}^{-1/2}$.
We also define the multiplication operator $S$ on $(L^2(\Omega ))^N$ by the matrix valued function
\begin{equation}
\label{esse}
S:= w^{-2}a^{-1/2} {\tilde a} a^{-1/2}\, .
\end{equation}
As it will be clear in the sequel, in order to compare $H$ and $\tilde H$ we shall also need the auxiliary operator $T^*ST$. Since
$$
\| S^{1/2}Tu \|^2_{L^2(\Omega ,\, g\, dx)}= \int_{\Omega }(  \tilde a \nabla u\cdot \nabla \bar u ) \tilde g    dx , \;\;\;\; u\in C_{\phi}[\cW],
$$
$T^*ST$ is the non-negative self-adjoint operator in $L^2(\Omega , g\, dx)$ canonically associated with the closure of the quadratic form
\[
\int_{\Omega }(  \tilde a \nabla u\cdot \nabla \bar u ) \tilde g    dx \; , \;\; u\in C_{\phi}[\cW] .
\]
So $\tilde H$ and $T^*ST$ have the same quadratic form, but they act on different Hilbert spaces:
$L^2(\Omega, \tilde g\, dx)$ and $L^2(\Omega , g\, dx)$ respectively.
It is easily seen that the operator $T^*ST$ is the pull-back to $\Omega$ via $\tilde\phi$ of the operator
\begin{equation}
\hat{L}:= \frac{\tilde g\circ \tilde\phi^{(-1)}}{g\circ \tilde\phi^{(-1)}} \tilde L.
\label{lhat}
\end{equation}
Thus we shall deal with the operators $L,\tilde L$ and $\hat L$ and the respective pull-backs $H$, $\tilde H$ and $T^*ST$. We shall repeatedly use the fact that these operators are pairwise unitarily equivalent.


Throughout this section  we assume that these operators have compact resolvent and that their eigenvalues satisfy the estimate
\begin{equation}
\lambda_n \geq C_1 n^{\frac{1}{\alpha}} \; , \qquad n\in\N,
\label{est}
\end{equation}
for some positive constants $\alpha$ and $C_1$.

\begin{remark}\label{remasym}
We recall that if $\Omega$ is a bounded domain with Lipschitz continuous boundary then (\ref{est}) is satisfied with  $\alpha =N/2$ (no restrictions on the boundary are required in the case of Dirichlet boundary conditions), see \cite{babula} for references.
\end{remark}

In the sequel we shall denote by $\lambda_n[E]$, $n\in\N$, the eigenvalues of a non-negative self-adjoint operator $E$ with compact resolvent, arranged in non-decreasing order and repeated according to multiplicity, and by $\psi_n[E]$,
$n\in\N$, a corresponding orthonormal sequence of eigenfunctions.

We introduce the following property which will be important in what follows.
\begin{definition}\label{growth}
Let $U$ be an open set in $\R^N$ and $\rho >0$ be a measurable function on $U$ and let $E$ be a non-negative self-adjoint operator on $L^2(U, \rho \, dx)$ with compact resolvent and $\dom(E)\subset W^{1,1}_{loc}(U)$. Let $q_0\in ]2,\infty ]$, $\gamma , C_2\in ]0,\infty [$. We say that $E$ satisfies property {\rm (P1)} with the parameters $q_0$, $\gamma$ and $C_2$ if
\[
\hspace{2.5cm} \| \psi_n[E]\|_{L^{q_0}(U , \rho \, dx)  }  \leq C_2\lambda_n[E]^{\gamma} \, , \;  \,  \hspace{3cm} {\rm (P1)}
\]
for all $n\in \N $ such that $\lambda_n[E]\ne 0$.
We say that $E$ satisfies property {\rm (P2)} with the parameters $q_0$, $\gamma$ and $C_2$ if
\[
\hspace{2.5cm} \| \nabla\psi_n[E]\|_{L^{q_0}(U , \rho \, dx)  }   \leq C_2\lambda_n[E]^{\frac{1}{2}+\gamma } \, , \;\; n\in\N \, . \hspace{2cm} {\rm (P2)}
\]
Finally, we say that $E$ satisfies property {\rm (P)} with the parameters $q_0$, $\gamma$ and $C_2$ if it satisfies both {\rm (P1)} and
{\rm (P2)} with these parameters.
\end{definition}

The next lemma involves the Schatten norms $\|\cdot\|_{\cC^r}$, $1\leq r\leq\infty$. For a compact operator $E$ on a Hilbert space they are defined by $\|E\|_{\cC^r}=(\sum_n \mu_n(E)^r)^{1/r}$, if $r<\infty$, and $\|E\|_{\cC^{\infty}}=\|E\|$,
where $\mu_n(E)$ are the singular values of $E$, {\it i.e.,} the non-zero eigenvalues of $(E^*E)^{1/2}$;
the Schatten space $\cC^r$, defined as the space of those compact operators for which the Schatten norm $\|\cdot\|_{\cC^r}$ is finite, is a Banach space; see Reed and Simon~\cite{RS} or Simon~\cite{S} for details.

Let $F:=TT^*$, $F_S:=S^{1/2}TT^*S^{1/2}$. It is well known that $\sigma(F)\setminus\{0\}=\sigma(H)\setminus\{0\}$ and similarly
$\sigma(F_S)\setminus\{0\}=\sigma(T^*ST)\setminus\{0\}$, see \cite[Theorem 2]{D}.
Moreover, we note that
\begin{equation}
\label{eich}
\tilde H= (\tilde a ^{1/2}\nabla ) ^{(*)_{\tilde g}}{\tilde a}^{1/2}\nabla = w^2({\tilde a}^{1/2}\nabla)^{*}w^{-2}{\tilde a}^{1/2}\nabla =w^2T^{*}ST\, ,
\end{equation}
and therefore the eigenvalues of the operator $wT^*STw$ coincide with the eigenvalues of  $\tilde H$.

\begin{lemma}\label{new}
$\ia$ Let $E$ be a non-negative self-adjoint operator on $L^2(\Omega , \rho\, dx)$ whose eigenvalues satisfy inequality (\ref{est}) for some $\alpha,C_1>0$. Assume that $E$ satisfies property {\rm (P1)} for some $q_0$, $\gamma$ and $C_2$. Then for large enough $k\in\N$, depending only on $\alpha$ and $\gamma$, there exists $c>0$ such that for all measurable functions $R$ on $\Omega$,
$$
 \|R(E+1)^{-k}\|_{\cC^2}\leq c  \|R\|_{L^{\frac{2q_0}{q_0-2}} (\Omega , \rho\, dx)  }.
$$
The constant $c$ depends only on $k$, $\alpha$, $\gamma$, $C_1$, $C_2$ and, if $\lambda_1[E]=0 $ and has multiplicity  $m$,  also on  $ \| \psi_i\|_{L^{q_0}(\Omega , \rho\, dx) }$, $i=1,\dots , m$.\newline
$\ib$  Assume that  $H$ (resp. $T^*ST$) satisfies property {\rm (P2)} for some $q_0$, $\gamma$ and $C_2$. Then for large enough $k\in\N$, depending only on  $\alpha$ and $\gamma$, there exists $c>0$ such that for all measurable matrix-valued functions $R$ on $\Omega$,
$$
\| R(F+1)^{-k}F^{1/2}\|_{\cC^2}\leq c  \| Ra^{1/2}\|_{L^{\frac{2q_0}{q_0-2}} (\Omega , g\, dx)  }.
$$
$$
( {\rm resp.} \quad \| R(F_S+1)^{-k}F_S^{1/2}\|_{\cC^2}\leq c  \| RS^{1/2}a^{1/2}\|_{L^{\frac{2q_0}{q_0-2}} (\Omega , g\, dx)  } .  \;\;)
$$
The constant $c$ depends only on $k$, $\alpha$, $\gamma$, $C_1$ and $C_2$.
\end{lemma}
{\em Proof.} We first prove statement (i). Assume for simplicity that $\lambda_1[E]\ne 0$.  We have
\begin{eqnarray}
\|R(E+1 )^{-k}\|_{\cC^2}^2&=& \sum_{n=1}^{\infty}\|R(E+1 )^{-k}\psi_n[E]\|^2_{L^2(\Omega , \rho \, dx)}\nonumber \\
&=& \sum_{n=1}^{\infty}(\lambda_n[E]+1)^{-2k}\|R\psi_n[E]\|^2_{L^2(\Omega , \rho \, dx)}\nonumber \\
&\leq& \|R\|_{L^{\frac{2q_0}{q_0-2}} (\Omega , \rho\, dx)  }   ^2\sum_{n=1}^{\infty}(\lambda_n[E]+1)^{-2k}\|\psi_n[E]\|^2_{L^{q_0}(\Omega , \rho \, dx)} \label{zero} \\
&\leq& c\|R\|_{L^{\frac{2q_0}{q_0-2}} (\Omega , \rho\, dx)  }  ^2\sum_{n=1}^{\infty}(\lambda_n[E]+1)^{-2k}\lambda_n[E]^{2\gamma} \nonumber \\
&\leq& c\|R\|_{L^{\frac{2q_0}{q_0-2}} (\Omega , \rho\, dx)  }  ^2 \, , \nonumber
\end{eqnarray}
provided $k$ is large enough. In case $\lambda_1[E]=0$ and has multiplicity $m$ one has simply to take into account the first $m$ summands in (\ref{zero}).

We now prove statement (ii). We only consider $F$, the operator $F_S$ is treated similarly. We note that $(\lambda_n[H]^{-1/2}T\psi_n[H])$ is an orthonormal basis of $\kernel (F)^{\perp}$. Hence
\begin{eqnarray*}
\| R(F +1)^{-k}F^{1/2}\| _{\cC^{2}}&= &\sum_{n=1}^{\infty}\lambda_n[E]^{-1}\|R(F+1 )^{-k}F^{1/2}T\psi_n[E]\|_{L^2(\Omega , g\, dx  )}   ^2\\
&=& \sum_{n=1}^{\infty}(\lambda_n[E]+1)^{-2k}\|RT\psi_n[E]\|_{L^2(\Omega , g\, dx  )}  ^2\\
&\leq&  \|Ra^{1/2}\|_{L^{\frac{2q_0}{q_0-2}} (\Omega , g\, dx)  }     ^2\sum_{n=1}^{\infty}(\lambda_n[E]+1)^{-2k}\|\nabla\psi_n[E]\|_{L^{q_0}(\Omega , g\, dx  )}  ^2 \\
&\leq&  c\|Ra^{1/2}\|_{L^{\frac{2q_0}{q_0-2}} (\Omega , g\, dx)  }    ^2\sum_{n=1}^{\infty}(\lambda_n[E]+1)^{-2k}\lambda_n[E]^{2\gamma+1} \\
&\leq & c\|Ra^{1/2}\|_{L^{\frac{2q_0}{q_0-2}} (\Omega , g\, dx)  }  ^2 \, ,
\end{eqnarray*}
provided $k$ is large enough. This completes the proof. $\hfill\Box$

\

Our stability estimates are expressed in terms of the following measure of vicinity of $\phi $ and $\tilde\phi$ (we recall that $w:=g^{1/2}{\tilde g}^{-1/2}$):
\begin{equation}\label{vicinity}
\delta _q(\phi ,\tilde\phi)= \delta_q^{(1)}(\phi ,\tilde\phi ) +\delta_q^{(2)}(\phi ,\tilde \phi),
\end{equation}
where
\begin{eqnarray*}
&& \delta _q^{(1)}(\phi ,\tilde \phi )=  \|w  -1\|_{L^{q}(\Omega ,  g\, dx)}  +\| w^{-1}-1 \|_{L^q(\Omega , g\, dx)}  \; ,\\
&& \delta _q^{(2)}(\phi ,\tilde \phi )= \| (S^{1/2}-S^{-1/2})a^{1/2}\|_{L^q(\Omega , g\, dx)} + \| (S-I)a^{1/2}\|_{L^q(\Omega , g\, dx)} .
\end{eqnarray*}

\begin{remark}\label{sobo}
Note that if we consider maps $\phi $ $\tilde \phi $ belonging to a family of transformations $\varphi $ satisfying the uniform estimate   $$ c^{-1}\le {\rm ess}\inf |{\rm det }\nabla \varphi |,\  \ \| \nabla \varphi\|_{L^{\infty }}\le c$$
 for a fixed $c>0$, and the coefficients $A_{ij}$ are Lipschitz continuous then
$$
\delta _q(\phi ,\tilde\phi) \le C \|\phi -\tilde\phi  \|_{W^{1,q}(\Omega )}.
$$
\end{remark}

\begin{theorem} {\bf (stability of resolvents)} Assume that the operators $H$ and $T^*ST$ satisfy properties {\rm (P1)} and {\rm (P2)}
and that $w^{-1}\tilde Hw$ satisfies property {\rm (P1)}, for the same  parameters $q_0$, $\gamma$ and $C_2$.
Then for all large enough $k\in\N$ depending only on $\alpha $ and $\gamma$ and for any $s>q_0(\alpha+2\gamma)/(q_0-2)$, there exists  $c>0$ such that
\begin{equation}
\|(w^{-1}\tilde Hw +1)^{-k} -(H+1 )^{-k}\|_{\cC^2}\leq c[1+\delta_s(\phi,\tilde\phi)]\delta _{\frac{2q_0}{q_0-2}}(\phi ,\tilde\phi ) .
\label{res2}
\end{equation}
The constant $c$ depends only on $\alpha, k,\gamma,s,C_1$, $C_2$ and, in the case of Neumann boundary conditions, also on $\|  g\|_{L^{q_0}(\Omega)}$.
\label{mainthm}
\end{theorem}
{\em Note.} The factor $1+\delta_s(\phi,\tilde\phi)$ appears for technical reasons and is not of importance for applications. \nl
{\em Proof.} We fix $k\in\N$ large enough so that part (i) of Lemma \ref{new} can be applied to the operators $H$, $T^*ST$ and $w^{-1}\tilde Hw $ and part (ii) of the same lemma can be applied to the operators $H$ and $T^*ST$. Since $w^{-1}\tilde Hw=wT^*STw$, we can write
\[
(w^{-1}\tilde Hw +1)^{-k} -(H+1)^{-k}=A+B,
\]
where
\[
A=(wT^*STw +1)^{-k} -(T^*ST+1 )^{-k} \; , \quad  B=(T^*ST +1)^{-k} -(T^*T+1 )^{-k}.
\]
We first estimate $A$ in terms of $\delta _q^{(1)}(\phi , \tilde \phi )$. We have
\begin{equation}
A=-  \sum_{i=0}^{k-1}(T^*ST+1)^{-i}  [ (T^*ST +1)^{-1} -(wT^*STw+1)^{-1}] (wT^*STw+1)^{-(k-1-i)}\, .
\label{a}
\end{equation}

First we estimate the terms in the sum (\ref{a}) corresponding to $i\leq [k/2]$. A direct computation shows that
\[
(T^*ST +1)^{-1} -(wT^*STw+1)^{-1} = D_1+D_2+D_3+D_4+D_5 \, ,
\]
where
\begin{eqnarray*}
D_1&=&(w-1)(wT^*STw+1 )^{-1} \\
D_2&=&(w-1)(wT^*STw+1 )^{-1}(w-1) \\
D_3&=&(wT^*STw+1 )^{-1} (w-1) \\
D_4&=&  (T^*ST+1 )^{-1}(w -w^{-1})(wT^*STw+1 )^{-1}(1-w) \\
D_5&=&(T^*ST+1 )^{-1}(w^{-1}-w)(wT^*STw+1 )^{-1} .
\end{eqnarray*}
Hence we need to estimate the terms $A_1,\ldots, A_5$ defined by
\[
A_j=\sum_{i=0}^{[k/2]}(T^*ST+1)^{-i} D_j (wT^*STw+1)^{-(k-1-i)}.
\]
Applying Lemma \ref{new} (i) for $E=wT^*STw$ we obtain that if $k$ is large enough, then
\begin{equation}
\label{porto}
 \|A_1\|_{\cC^2}\leq c  \|w-1\|_{L^{\frac{2q_0}{q_0-2}}(\Omega , g\, dx)}.
\end{equation}
Now,  applying \cite[Lemma 4.5]{babula} with $p=q_0/(q_0-2)$ we get
$\|(w-1)(wT^*STw+1 )^{-1}\| \leq c\|w-1\|_{L^s(\Omega , g\, dx)}$, for any $s>q_0(\alpha+2\gamma)/(q_0-2)$, hence
$\|A_2\|_{\cC^2}\leq c  \|w-1\|_{L^s(\Omega , g\, dx)}\|w-1\|_{ L^{2q_0/(q_0-2)}(\Omega ,g\, dx )  }$. The remaining terms $A_3$, $A_4$ and $A_5$ are estimated similarly.

In order to estimate the terms in the sum (\ref{a}) corresponding to $i> [k/2]$ it is possible to proceed as above by swapping $(T^*ST +1)$ and
$(wT^*STw+1)$ and using  the following
decomposition
\[
(T^*ST +1)^{-1} -(wT^*STw+1)^{-1} = D_1'+D_2'+D_3'+D_4'+D_5' \, ,
\]
where
\begin{eqnarray*}
D_1'&=&(1-w^{-1})(T^*ST+1 )^{-1} \\
D_2'&=&(w^{-1}-1)(T^*ST+1 )^{-1}(1-w^{-1}) \\
D_3'&=&(T^*ST+1 )^{-1} (1-w^{-1}) \\
D_4'&=&  (w^{-1}-1)(T^*ST+1 )^{-1}(w^{-1}-w)(wT^*STw+1 )^{-1} \\
D_5'&=&(T^*ST+1 )^{-1}(w^{-1}-w)(wT^*STw+1 )^{-1} .
\end{eqnarray*}

We now consider the term $B$. We write
\[
B=\sum_{i=0}^{k-1}(T^*ST+1)^{-i}  [ (T^*ST +1)^{-1} -(T^*T+1)^{-1}] (T^*T+1)^{-(k-1-i)} \, .\\
\]
Let $B_i$ denote the $i$th summond. We have \cite{babula}
\[
(T^*ST +1)^{-1} -(T^*T+1)^{-1} =T^*S^{1/2}(F_S+1)^{-1}(S^{-1/2}-S^{1/2})(F+1)^{-1}T \, .
\]
It is also known \cite{D} that
$T(T^*T+1)^{-m}=(TT^*+1)^{-m}T$, $m\in\N$, with a similar relation, of course, for $S^{1/2}T$. Hence
\begin{eqnarray*}
B_i&=& (T^*ST+1)^{-i}T^*S^{1/2}(F_S+1)^{-1}(S^{-1/2}-S^{1/2})(F+1)^{-1}T(T^*T+1)^{-(k-1-i)}\\
&=&T^*S^{1/2} (F_S+1)^{-i-1} (S^{-1/2}-S^{1/2})(F+1)^{-k-i}T.
\end{eqnarray*}
Using polar decomposition for $S^{1/2}T$ we note that $\|T^*S^{1/2}(F_S+1)^{-i-1}\|\le 1$. Using also polar decomposition for $T$
and applying Lemma \ref{new} (ii) for $F$ we therefore obtain that for  $i\leq [k/2]$ there holds
\begin{eqnarray*}
\|B_i\|_{\cC^2} &\leq&\|T^*S^{1/2}(F_S+1)^{-i-1}\| \|(S^{-1/2}-S^{1/2})(F+1)^{-(k-i)}F^{1/2}\|_{\cC^2}\\
&\leq&\|(S^{-1/2}-S^{1/2})(F+1)^{-(k-i)}F^{1/2}\|_{\cC^2} \\
&\leq& c\|(S^{-1/2}-S^{1/2})a^{1/2}\|_{L^{\frac{2q_0}{q_0-2}}(\Omega , g\, dx)},
\end{eqnarray*}
provided $k$ is large enough. For the terms with $i> [k/2]$, we argue similarly (but now use $F_S$ instead of $F$) and obtain
\[
\|B_i\|_{\cC^2} \leq c\|(S-I)a^{1/2}\|_{L^{\frac{2q_0}{q_0-2}}(\Omega , g\, dx)} .
\]
This concludes the proof of the theorem.  $\hfill\Box$\\

As in \cite{babula}, from Theorem~\ref{mainthm} we immediately deduce the following  theorem.
\begin{theorem} {\bf (stability of eigenvalues)} Assume that the operators $H$ and $T^*ST$ satisfy properties {\rm (P1)} and {\rm (P2)}
and that $w^{-1}\tilde Hw$ satisfies property {\rm (P1)}, for the same  parameters $q_0$, $\gamma$ and $C_2$.
Then for all large enough $k\in\N$ depending only on $\alpha, \gamma$ and for any $s>q_0(\alpha+2\gamma)/(q_0-2)$
there exists  $c>0$ such that
\begin{equation}
\left( \sum_{n=1}^{\infty}\left[
(\lambda_n[\tilde L] +1)^{-k} -  (\lambda_n[L] +1)^{-k}\right]^2\right)^{1/2}
\leq c [1+\delta_s(\phi,\tilde\phi)]\delta_{\frac{2q_0}{q_0-2}}(\phi ,\tilde\phi).
\label{res2'}
\end{equation}
The constant $c$ depends only on $\alpha, k,\gamma,s,C_1, C_2$ and, in the case of Neumann boundary conditions, also on $\|  g\|_{L^{q_0}(\Omega)}$.
\label{thm:series}
\end{theorem}

In order to estimate the variation of the eigenfunctions, we need the following

\begin{lemma}\label{boundedop}
Let $A,B$ be compact self-adjoint and positive operators in a Hilbert space ${\mathcal {H}}$. Let $\lambda_n$, $\mu_n$, $n\in {\N}$ be the eigenvalues of $A, B$ respectively. Let $\phi_n, \psi_n$, $n\in {\N}$ be orthonormal sequences of eigenfunctions corresponding to $\lambda_n, \mu_n$ respectively. Let $\nu $ be an eigenvalue of $A$, $\Lambda=\{n\in {\N}:\ \lambda_n=\nu \}$ and $d=\min \{|\lambda_i-\nu |:\, i\in {\N}\setminus \Lambda \}$. Let $P, Q$ be the orthogonal projectors of ${\mathcal{H}}$ onto ${\rm span }\{\phi_n:\ n\in \Lambda \}$ and ${\rm span }\{\psi_n:\, n\in \Lambda \} $, respectively.

If $\| A-B\|< d/2$ then $\| P-Q\| < \frac{2(1+|\Lambda|)}{d}\| A-B \|$.
\end{lemma}
{\it Proof.} Note that by the min-max Principle it follows that $|\lambda_i -\mu_i|\le \| A-B \|$ for all $i \in {\N}$; thus, if  $\| A-B\|< d/2$ then $|\mu_n-\nu |< d/2 $ for all $n\in \Lambda$ and $\|\mu_i -\nu \|> d/2$ for all $i\in {\N}\setminus \Lambda$.

Let $u\in {\rm Ran}(P)$, $\| u\| \le 1$. Then
\begin{eqnarray*}
\| A-B\|^2 & \geq  & \| Au-Bu\|^2=\| \nu u-B u \| ^2\\
&= & \| \nu \sum_{i=1}^{\infty }\inprod{u}{\psi_i }\psi_i- \sum_{i=1}^{\infty }\mu_i\inprod{u}{\psi_i }\psi_i  \|^2\\
&\geq  & \sum_{i\notin \Lambda}(\nu -\mu_i)^2|\inprod{u}{\psi_i}|^2> \frac{d^2}{4}  \sum_{i\notin \Lambda}|\inprod{u}{\psi_i}|^2 ,
\end{eqnarray*}
that is $\| (I-Q)u\| <\frac{2}{d}\| A-B\| $. Thus
\begin{equation}\label{pr1}
\|(I-Q )P \| <\frac{2}{d}\| A-B\|\, .
\end{equation}
Now, let $n\in \Lambda$. Then
\begin{eqnarray*}
\| A -B\|^2  & \geq  & \| A\psi_n -B\psi_n\|^2 =\| \sum_{i=1}^{\infty } \lambda_i\inprod{\psi_n}{\phi_i}\phi_i-\mu_n\psi_n \|^2\\
& =&   \sum_{i=1}^{\infty } (\lambda_i-\mu_n)^2|\inprod{\psi_n}{\phi_i}|^2 \geq \sum_{i\notin \Lambda} (\lambda_i-\mu_n)^2|\inprod{\psi_n}{\phi_i}|^2 \\
&> & \frac{d^2}{4}\sum_{i\notin \Lambda} |\inprod{\psi_n}{\phi_i}|^2=\frac{d^2}{4}\| (I-P)\psi_n \|^2.
\end{eqnarray*}
Hence
\begin{equation}\label{pr2}
\| Q(I-P) \| = \| (I-P)Q \| \le \frac{2|\Lambda|\| A-B\|}{d} .
\end{equation}
The proof follows by combining (\ref{pr1}) and (\ref{pr2}). \hfill $\Box$\\

In the following theorem it is understood that $\psi_k[L]$ and $\psi_k[\tilde L]$ are extended by zero outside $\phi(\Omega)$ and $\tilde\phi(\Omega)$ respectively.
\begin{theorem} {\bf (stability of eigenfunctions)}
\label{thm:eigen}
Assume that the operators $H$ and $T^*ST$ satisfy properties {\rm (P1)} and {\rm (P2)}
and that $w^{-1}\tilde Hw$ satisfies property {\rm (P1)}, for the same  parameters $q_0$, $\gamma$ and $C_2$.
Let $\lambda$ be an eigenvalue of $L$ (resp. $\tilde L$) of multiplicity $m$ and let $n\in\N$ be such that
$\lambda=\lambda_n=\ldots =\lambda_{n+m-1}$.
Then for any $s>q_0(\alpha +2\gamma )/(q_0-2) $  there exists  $c>0$ depending only on $\alpha , q_0, \gamma , C_1,C_2, \lambda_{n-1} , \lambda $, $ \lambda_{n+m}$ and, in case of Neumann boundary conditions, $\| g\|_{L^{q_0}(\Omega )}$ such that the following is true: if
$
[1+\delta_s(\phi,\tilde\phi)]\delta_{2q_0/(q_0-2)}(\phi ,\tilde\phi)\leq c^{-1}
$
and $\psi_n[\tilde L], \dots $, $\psi_{n+m-1}[\tilde L]$ (resp. $\psi_n[L], \dots $, $\psi_{n+m-1}[L]$) are orthonormal  eigenfunctions of $\tilde L$  in $L^2(\tilde\phi (\Omega ))$ (resp. $L$  in $L^2(\phi(\Omega) )$), then there exist orthonormal eigenfunctions
$ \psi_n[L], \dots , \psi_{n+m-1}[L]$ of $L$ in $L^2(\phi(\Omega))$  (resp. $ \psi_n[\tilde L], \dots , \psi_{n+m-1}[\tilde L]$ of $\tilde L$ in $L^2(\tilde\phi (\Omega))$ ) such that
\begin{eqnarray}\label{res2''}
\label{holder2}\lefteqn{
 \| \psi_{l}[L]-\psi_{l}[\tilde L]  \|_{L^2(\phi(\Omega) \cup \tilde\phi (\Omega ))} \le  c \big(
[1+\delta_s(\phi,\tilde\phi)] \delta_{\frac{2q_0}{q_0-2} }(\phi ,\tilde\phi) +}\nonumber \\
& &
+  \| g^{1/2}\psi_l[L]\circ\phi -\tilde g^{1/2}\psi_l[L]\circ \tilde \phi  \|_{L^{2 }(\Omega ) }+
  \| g^{1/2} \psi_l[ \tilde L]\circ\phi -\tilde g^{1/2}\psi_l[\tilde L]\circ  \tilde\phi  \|_{L^{2 }(\Omega ) }  \big),\nonumber \\
 \end{eqnarray}
 for all $l=n, \dots , n+m-1$.
\end{theorem}
{\it Proof.} Let $k\in {\N}$ be large enough so that estimate (\ref{res2}) holds.
Let $\lambda =\lambda_n=\dots =\lambda_{n+m-1}$ be an eigenvalue of $L$ of multiplicity $m$. By applying Lemma~\ref{boundedop} with
${\mathcal{H}}=L^2(\Omega , g\, dx)$, $A=(H+1)^{-k}$, $B=(w^{-1}\tilde Hw +1)^{-k}$ and $\nu = (\lambda +1)^{-k}$, it follows that there exists $c>0$ as in the statement such that if   $[1+\delta_s(\phi,\tilde\phi)]\delta_{\frac{2q_0}{q_0-2}}(\phi ,\tilde\phi)\le c^{-1}$ then
\begin{equation}
\label{estproj}
\| P-Q\| \le c [1+\delta_s(\phi,\tilde\phi)]\delta_{\frac{2q_0}{q_0-2}}(\phi ,\tilde\phi),
\end{equation}
where $P,Q$ are the orthononal projectors in $L^2(\Omega , g\, dx)$ as in Lemma~{\ref{boundedop}}.

Now, given eigenfunctions $\psi _l [\tilde L]$ as in the statement, we set $\psi_l[\tilde H]=\psi_l[\tilde L]\circ \tilde \phi $ and we note that
$w^{-1}\psi_l[\tilde H]$ are eigenfunctions of $w^{-1}\tilde H w$.
Proceeding as in the proof of \cite[Thm.~5.6]{babula}, using the Selection Lemma \cite[Lemma.~5.4]{babula} and estimate (\ref{estproj}), we have  that by possibly enlarging $c$,  if $[1+\delta_s(\phi,\tilde\phi)]\delta_{2q_0/(q_0-2)}(\phi ,\tilde \phi )<c^{-1}$   there exist eigenfunctions  $\psi_n[H]$, $\dots ,$ $ \psi_{n+m-1}[H]$  such that
\begin{equation}\label{autof0}
\| \psi_l[H] -w^{-1}\psi_l[\tilde H]\|_{L^2(\Omega , g\, dx)}\le c  [1+\delta_s(\phi,\tilde\phi)]\delta_{2q_0/(q_0-2)}(\phi ,\tilde \phi ),
\end{equation}
for all $l=n, \dots , n+m-1$.
We set $\psi _l[L]=\psi _l[H]\circ \phi^{(-1)} $ for all $l=n, \dots , n+m-1$.
By changing variables in the left-hand side of (\ref{autof0}) we obtain
$$
\| \psi_l[L] -w^{-1}\circ \phi^{(-1)}\psi_l[\tilde H]\circ \phi^{(-1)}\|_{L^2(\phi (\Omega ))}\le c  [1+\delta_s(\phi,\tilde\phi)]\delta_{2q_0/(q_0-2)}(\phi ,\tilde \phi ),
$$
hence
\begin{eqnarray*}
\| \psi_l[L]-\psi_l[\tilde L]\| _{L^2(\phi (\Omega ))}\le c  [1+\delta_s(\phi,\tilde\phi)]\delta_{2q_0/(q_0-2)}(\phi ,\tilde \phi )\\
+\| \psi_l[\tilde L] -w^{-1}\circ \phi^{(-1)}\psi_l[\tilde H]\circ \phi^{(-1)}\|_{L^2(\phi (\Omega ))}
\end{eqnarray*}
and similarly
\begin{eqnarray*}
\| \psi_l[\tilde L]-\psi_l[ L]\| _{L^2(\tilde \phi (\Omega ))}\le c  [1+\delta_s(\phi,\tilde\phi)]\delta_{2q_0/(q_0-2)}(\phi ,\tilde \phi )\\
+\| \psi_l[ L] -w\circ \tilde \phi^{(-1)}\psi_l[ H]\circ\tilde\phi^{(-1)}\|_{L^2(\tilde \phi (\Omega ))}\, ,
\end{eqnarray*}
for all $l=n, \dots , n+m-1$. Estimate (\ref{holder2}) follows by combining the previous two inequalities and changing variables in integrals again.
 If $\lambda$ is an eigenvalue of $\tilde L$ we work similarly.\hfill $\Box$

\section{Applications}
\label{applicationsec}

In this section we apply Theorems \ref{thm:series} and \ref{thm:eigen} in order to obtain explicit stability estimates in terms of Lebesgue measure.
This will be carried out by showing that condition (P) is satisfied in suitable classes of domains and by constructing appropriate transformations
$\phi$.

\subsection{Spectral stability for smooth and Lipschitz domains}

In this subsection we consider bounded domains  $\Omega$ in $\R^N$ of class $C^{m,1}$ for  $m=0,1$, {\it i.e.,} bounded domains which are locally subgraphs of $C^{m,1}$ functions. In this context, domains of class $C^{1,1}$ represent the smooth case.

\begin{theorem}\label{regex}The following statements hold.
\begin{itemize}\item[$\ia$] Let $\Omega $ be a  bounded domain in $\R^N$  of class $C^{1,1}$ and
let $A_{ij}$ be Lipschitz  functions defined on $\Omega $ satisfying (\ref{ellip}). Then the operator (\ref{classic})
subject either to Dirichlet or Neumann boundary conditions on $\Omega $ satisfies property (P) with $q_0=\infty $ and $\gamma =N/4 $.
\item[$\ib$] Let $\Omega $ be a bounded domain in $\R^N$  of class $C^{0,1}$ and
let $A_{ij}$ be measurable  functions defined on $\Omega $ satisfying (\ref{ellip}). Then the operator (\ref{classic})
subject either to Dirichlet or Neumann boundary conditions on $\Omega $ satisfies property (P) with some $q_0>2 $ and $\gamma =N(q_0-2)/(4q_0) $.
\item[$\ic$] The Laplace operator subject to Dirichlet boundary conditions on a bounded domain  in ${\mathbb{R}}^N$ of class
$C^{0,1}$ satisfies property (P) with some $q_0>4$ if $N=2$ and some $q_0>3$ if $N\geq 3$.
\end{itemize}
\end{theorem}

Statement $\ia$ is well-known, see \cite{babula} for references. For a proof of statement $\ib$ we refer, {\it e.g.}, to \cite[Remark~6.5]{babula}.
Statement $\ic$ is a consequence of Jerison and Kenig~\cite{jeke}.

\begin{definition}
Let $V$ be a bounded open cylinder, i.e., there exists a rotation  $R$  such that $R(V)=W\times ]a,b[$, where $W$ is a bounded convex open set in $\R^{N-1}$. Let $M,\rho>0$. We say that a bounded domain $\Omega\subset\R^N$ belongs to $\cC^{m,1}_M(V,R,\rho)$ if
$\Omega$ is of class $C^{m,1}$  and there exists a function $g\in C^{m,1}(\overline{W})$ such that $a+\rho\leq g\le b$,   $| g |_{m,1}:=\sum_{0<|\alpha |\le m+1}\| D^{\alpha }g \|_{L^{\infty }(W )}\le M$,
 and
\begin{equation}
R(\Omega\cap V)=\{ (\bar{x},x_N) \; : \;  \bar{x}\in W \, ,\,  a<x_N<g(\bar{x}) \}.
\label{C11}
\end{equation}
\label{def:C11}
\end{definition}
In the following theorem we denote by $\lambda_n[L]$, $\lambda_n[\tilde L]$  the eigenvalues of the operator (\ref{classic}) subject to Dirichlet
or Neumann boundary conditions on $\Omega $ and $\tilde \Omega $ respectively.  Similarly, the eigenfunctions are denoted by $\psi_n[L]$ and $\psi_n[\tilde L]$. Moreover, by $V_{\rho }$ we denote the set $\{x\in V:\ d(x, \partial V )>\rho \}$.

\begin{theorem} Let $A_{ij}$, $i,j=1,\dots ,N$ be measurable functions defined on $\R^N$ satisfying $A_{ij}=A_{ji}$ and the ellipticity condition
(\ref{ellip}). Let $\Omega \in \cC^{0,1}_M(V,R,\rho)$. Then there exists $2<q_0\leq\infty$ such that the following statements hold:
\begin{itemize}
\item[$\ia$]
For all large enough $k\in\N $ there exists $c_1>0$ such that
\begin{equation}
\left(\sum_{n=1}^{\infty}\left| (
{\lambda}_n[L] +1)^{-k} - ({\lambda}_n[\tilde L] +1)^{-k}  \right|^2\right)^{1/2}
\leq c_1 |\Omega \vartriangle \tilde \Omega  |^{\frac{1}{2}-\frac{1}{q_0}} ,
\end{equation}
for all $\tilde \Omega \in  \cC^{0,1}_M(V,R,\rho) $ such that $\tilde\Omega \cap (V_{\rho})^{c}=\Omega \cap (V_{\rho})^{c}$.
\item[$\ib$] Let $\lambda [L] $ be an eigenvalue of multiplicity $m$ and let $n\in \N$ be such that
$\lambda [L ] =\lambda_n[L ]=\dots =\lambda_{n+m-1}[L]$.
There exists $c_2>0$ such that the following is true: if
$\tilde \Omega \in  \cC^{0,1}_M(V,R,\rho) $, $\Omega \cap (V_{\rho})^{c}=\tilde \Omega \cap (V_{\rho})^{c}$,
$
|\Omega \vartriangle \tilde \Omega |\le c_2^{-1},
$
then, given orthonormal eigenfunctions
$\psi_n[\tilde L ], \dots , $ $ \psi_{n+m-1}[\tilde L]$ in $L^2(\tilde \Omega )$, there exist corresponding  orthonormal eigenfunctions $ \psi_n[L ],$ $ \dots ,$   $  \psi_{n+m-1}[L ] $ in $L^2(\Omega  )$
such that
$$
\| \psi_n[L ]-\psi_n[\tilde L ] \|_{L^2(\Omega \cup\tilde \Omega)}\le c_2|\Omega \vartriangle\tilde \Omega |^{\frac{1}{2}-\frac{1}{q_0}  }.
$$
\end{itemize}
If in addition $A_{ij}\in C^{0,1}(\R^N)$ and $\Omega, \tilde \Omega \in  \cC^{1,1}_M(V,R,\rho)$ then statements (i) and (ii) hold with $q_0=\infty $.
\label{thm:C11}
\end{theorem}
{\em Proof.}  Let $\tilde \Omega  \in \cC^{0,1}_M(V,R,\rho )$. By \cite[Lemma~7.4]{babula} there exists a bi-Lipschitz map $\Phi $ from  $\Omega $ onto $\tilde \Omega $ such that
\begin{equation}\label{tauphi}\tau^{-1} \le {\rm ess}\inf_{\Omega }|{\rm det }\nabla \Phi | \quad  {\rm and }\quad \| \nabla \Phi \| _{L^{\infty }(\Omega )}\le \tau\, ,
\end{equation}
where $\tau >0$ depends only on $N,V,M,\rho $ and such that there exists $\hat \Omega \subset \Omega $ satisfying the following properties:
\begin{equation}\label{omegahat}\Phi (x) =x\ {\rm for}\ {\rm all}\ x\in \hat \Omega \ {\rm and}\ |\Omega \setminus \hat\Omega |\le 2 |\Omega \vartriangle \tilde\Omega |\, .
\end{equation}
As in \cite[Theorem~7.3]{babula} we apply Theorems~\ref{thm:series} and \ref{thm:eigen} with $\phi =Id$ and $\tilde\phi =\Phi$. By Remark~\ref{remasym}  condition (\ref{est}) is satisfied for  $\alpha =N/2$. Moreover, by Theorem~\ref{regex} and \cite{babula} the operators $L$, $\tilde L$ and $T^*ST$ satisfy property (P) for some $q_0>2$, hence also $H, \tilde H $ satisfy property (P) and $w^{-1}\tilde H w$ satisfies property (P1). Thus Theorems~\ref{thm:series} and \ref{thm:eigen} apply  and estimates (\ref{res2'}), (\ref{res2''}) hold. By (\ref{tauphi}) and (\ref{omegahat}) it follows that
$$
\delta_{2q_0/(q_0-2)}(\phi ,\tilde\phi )\le c|\Omega\vartriangle \tilde\Omega |^{\frac{1}{2}-\frac{1}{q_0}}
$$
which combined with (\ref{res2'}), (\ref{res2''}) easily implies the validity of statements $\ia$ and $\ib$ (the last two terms in the right-hand side of  (\ref{holder2}) are estimated by means of the H\"{o}lder inequality and property (P1)). In the case of open sets of class
$C^{1,1}$ it is enough to observe that by Theorem~\ref{regex} and \cite{babula} it is possible to choose $q_0=\infty $ and proceed as above.
 \hfill $\Box$

\subsection{An abstract regularity theorem}\label{sub1}

We prove a theorem on the regularity of eigenfunctions of a general operator $H$ which will be used in the next subsection. This theorem is a  generalization of \cite[Thm.~5.1]{bulaneu} which was concerned with domains satisfying a uniform cone condition. The theorem has two main assumptions: a general multiplicative Sobolev inequality (which is an assumption on the underlying domain $\Omega$ and replaces the standard multiplicative Sobolev inequality used in \cite{bulaneu}) and an a priori estimate on the operator $H$. More precisely, we need to consider the following properties:\\

{\bf (A) Sobolev inequality.} Let $m\in\N$, $M>0$. If $1\leq p,q\leq\infty$ and $\beta$ is a multi-index of length
$|\beta|<m$ such that
\[
\frac{1}{q}\geq\frac{1}{p} -\frac{m-|\beta|}{M} ,
\]
where, in the case of equality, $1<p<q<\infty$,
then there exists $\tau=\tau(m,\beta,M,p,q)$ $\in ]0,1]$ and $C_4=C_4(m,\beta,M,p,q,\Omega)$ such that for all $u\in W^{m,p}(\Omega)$,
\begin{equation}
\| D^{\beta}u \|_{L^q(\Omega)} \leq C_4 \| u\|_{W^{m,p}(\Omega)}^{\tau} \|u\|_{L^p(\Omega)}^{1-\tau}.
\label{sob_cusp}
\end{equation}

{\bf (B) A priori estimate.} Let $m\in\N$, $1< p_0<\infty$ and $H:\dom(H)\to L^1_{loc}(\Omega)$ where $\dom(H)\subset L^{p_0}(\Omega)$. For all $p_0\leq p<\infty$ there exists $A_p<\infty$ such that if $u\in\dom(H)$ and $Hu\in L^p(\Omega)$ then $u\in W^{m,p}(\Omega)$ and
\begin{equation}
\|u\|_{W^{m,p}(\Omega)}\leq A_p \|Hu\|_{L^p(\Omega)} .
\label{apriori}
\end{equation}
Then, following the bootstrap argument used in  \cite{bulaneu}, we prove the following
\begin{theorem}
Let $1<p_0<\infty$ and let $H$ be an operator with $\dom(H)\subset L^{p_0}(\Omega)$. Assume that the Sobolev inequality (A) and the a priori estimate (B) are satisfied for some
$m\in\N$, $M>0$ . Assume further that $\tau(m,0,M,p,q)=\frac{M}{m}(\frac{1}{p}-\frac{1}{q})$.
Then  for any eigenfunction $\phi$ of $H$, $H\phi=\lambda\phi$, the following statements hold:
\begin{eqnarray}
\ia&& \mbox{For any $p_0\leq p<\infty$  $\phi \in W^{m,p}(\Omega )$ and there exists $B_p<\infty$ such that}\nonumber \\
&&\qquad\|\phi\|_{W^{m,p}(\Omega)}\leq B_p |\lambda|^{1+\frac{M}{m}(\frac{1}{p_0}-\frac{1}{p})}  \|\phi\|_{L^{p_0}(\Omega)}.\;\label{reg1} \\[0.2cm]
\ib&& \mbox{Let $\beta$ be a multi-index with $|\beta|<m$ and define}\nonumber\\
&& \quad\rho=\inf\Big\{ \tau(m,\beta,M,p,\infty)+\frac{M}{m}\Big(\frac{1}{p_0}-\frac{1}{p}\Big)   \; : \; p>\max\Big\{\frac{M}{m-|\beta|},p_0\Big\}\Big\} .\nonumber\\
&&\mbox{Then for any $\eta>0$ there exists $B_{\mu,\eta}<\infty$ such that }\nonumber\\
&&\qquad\|D^{\beta}\phi\|_{L^{\infty}(\Omega)}\leq B_{\mu,\eta} (1+|\lambda|)^{\rho+\eta}\|\phi\|_{L^{p_0}(\Omega)}.
\label{reg2}
\end{eqnarray}
\label{bootstrap_thm}
\end{theorem}
\begin{remark}\label{remboot} {\rm It is immediate that if
\[
\tau(m,\beta,M,p,\infty)=\frac{A+\alpha M/p}{A +\alpha(m-|\beta|)}\, ,
\]
for some $A,\alpha>0$, and $|\beta |<m$ and $p_0\leq M/(m-|\beta|)$, then by (\ref{reg2}) it follows that
\begin{equation}
\|D^{\beta } \phi\|_{L^{\infty}(\Omega)}\leq B_{|\beta |,\eta} (1+|\lambda|)^{\frac{1}{m}(\frac{M}{p_0}+|\beta|)+\eta}\|\phi\|_{L^{p_0}(\Omega)} ,
\label{reg3}
\end{equation}
for any $\eta>0$.}
\end{remark}
{\em Proof.} We set $s(q)=M q/(M -mq)$ if $0\leq q< M /m$, and $s(q)=\infty$ if $M /m\leq q\leq\infty$. Since
$$\lim_{k\to \infty}\underbrace{s (\dots(s (q))\dots )}_{k}=\infty $$ one can apply the bootstrap argument used in \cite[Theorem 5.1]{bulaneu} and prove that $\phi\in L^p(\Omega)$ for any $p_0\leq p\leq\infty$; see \cite[Remark 5.9]{bulaneu}.
Applying (\ref{apriori}) we obtain
\begin{equation}
\|\phi\|_{W^{m,p}(\Omega)}\leq A_p |\lambda|\|\phi\|_{L^p(\Omega)}  \; , \quad p_0\leq p<\infty\, .
\label{sim}
\end{equation}

If $p=p_0$, then (\ref{reg1}) is an immediate consequence of (\ref{sim}), so we assume that $p_0<p<\infty$.
Let us define $\sigma  (t)=M t/(M+m t)$, for all $t\geq 0$. Note that   $\sigma (t)=s^{(-1)}(t)$ for all $t\geq 0$. We define a sequence $(p_k)_{k\geq 1}$ by
\[
p_1=p\; , \qquad p_{k+1}=\max\{ p_0,\frac{1}{2}(\sigma(p_k)+p_k)\}\; .
\]
We note that
\[
\frac{1}{2}(\sigma(p_k)+p_k) =\Big( M +\frac{mp_k}{2}\Big)(M +mp_k)^{-1} p_k\leq\nu p_k\, ,
\]
where $\nu=(M +\frac{mp_0}{2})(M+mp_0)^{-1}<1$; hence $p_k=p_0$ from some $k$ onwards. Let $\kappa$ be the first such $k$.
We then have $p=p_1>p_2>\ldots  >p_{\kappa}=p_0$. Moreover, $p_{k+1}>\sigma(p_k)$, $k=1,\ldots,\kappa-1$. Inverting we then obtain $p_k<s(p_{k+1})$, $k=1,\ldots,\kappa-1$. Applying (\ref{sob_cusp}) for $q=p_k$, $p=p_{k+1}$, $\beta=0$ and (\ref{sim}) we obtain, with
$\tau_k =\tau(m,0,M,p_{k+1},p_k)$ and $c_1(k)=C_4(m,0,M,p_{k+1},p_k)A_{p_{k+1}}^{\tau_k}$,
\begin{eqnarray*}
\|\phi\|_{L^{p_k}(\Omega)} &\leq& C_4(m,0,M,p_{k+1},p_k) \|\phi\|_{W^{m,p_{k+1}}(\Omega)}^{\tau_k} \|\phi\|_{L^{p_{k+1}}(\Omega)}^{1-\tau_k} \\
&\leq& C_4(m,0,M,p_{k+1},p_k) A_{p_{k+1}}^{\tau_k}
|\lambda|^{\tau_k}\|\phi\|_{L^{p_{k+1}}(\Omega)}^{\tau_k} \|\phi\|_{L^{p_{k+1}}(\Omega)}^{1-\tau_k} \\
&=& c_1(k)|\lambda|^{\tau_k}\|\phi\|_{L^{p_{k+1}}(\Omega)} .
\end{eqnarray*}
Hence
\begin{eqnarray*}
\|\phi\|_{L^{p_1}(\Omega)} &\leq& c_2(1)|\lambda|^{\tau_1}\|\phi\|_{L^{p_{2}}(\Omega)} \\
&\leq& c_1(1)c_1(2)|\lambda|^{\tau_1+\tau_2}\|\phi\|_{L^{p_{3}}(\Omega)},
\end{eqnarray*}
and iterating
\[
\|\phi\|_{L^{p_1}(\Omega)} \leq c_2(p)|\lambda|^{\tau_1+\tau_2+\cdots +\tau_{\kappa -1}}\|\phi\|_{L^{p_{\kappa}}(\Omega)} ,
\]
where $c_2(p)=\prod_{i=1}^{\kappa-1}c_1(i)$. Recalling that $p_1=p$, $p_{\kappa}=p_0$ and $\tau_1+\ldots+\tau_{\kappa-1}=
\frac{M }{m}(\frac{1}{p_{0}}-\frac{1}{p})$, this takes the form
\begin{equation}
\|\phi\|_{L^{p}(\Omega)} \leq c_2(p)|\lambda|^{\frac{M }{m}(\frac{1}{p_{0}}-\frac{1}{p})}\|\phi\|_{L^{p_0}(\Omega)} .
\label{sim1}
\end{equation}
Plugging this back to (\ref{sim}) we obtain (\ref{reg1}), with $B_p=c_2(p)A_p\, $.

We now prove (ii). Let $|\beta|<m$, $p>\max\big\{\frac{M}{m-|\beta|},p_0\big\}$ and $\tau=\tau(m,\beta,M,p,\infty)$, $C_4=C_4(m,\beta,M,p,\infty ,\Omega)$. By (\ref{sob_cusp}), (\ref{reg1}) and (\ref{sim1}) we then have
\begin{eqnarray}
\| D^{\beta}\phi \|_{L^{\infty}(\Omega)} &\leq&  C_4\| \phi\|_{W^{m,p}(\Omega)}^{\tau} \|\phi\|_{L^p(\Omega)}^{1-\tau}\nonumber\\
&\leq& C_4\Big( B_p |\lambda|^{1+\frac{M }{m}(\frac{1}{p_0}-\frac{1}{p})}  \|\phi\|_{L^{p_0}(\Omega)} )\Big)^{\tau}\times \nonumber\\
&& \qquad \Big(c_2(p)|\lambda|^{\frac{M}{m}(\frac{1}{p_{0}}-\frac{1}{p})}\|\phi\|_{L^{p_0}(\Omega)} \Big)^{1-\tau}\nonumber\\
&= &c_3(p)|\lambda|^{\rho} \|\phi\|_{L^{p_0}(\Omega)} ,
\label{arl}
\end{eqnarray}
where $c_3(p)=C_4B_p^{\tau}c_2(p)^{1-\tau}$ and
\[
\rho=\tau\left(1+\frac{M }{m}\left(\frac{1}{p_0}-\frac{1}{p}\right)\right) +(1-\tau)\frac{M }{m}\left(\frac{1}{p_{0}}-\frac{1}{p}\right)
=\tau+\frac{M }{m}\left(\frac{1}{p_0}-\frac{1}{p}\right).
\]
Optimizing over $p$ we obtain (\ref{reg2}).$\hfill\Box$

\subsection{Spectral stability for domains with outward cusps }\label{applsec}

Let $0<\alpha<1$. Let $\Omega\subset\R^N$, $N\geq 2$, be a domain the boundary of which
is $C^2$ apart from a single outward cusp. More precisely we assume that
\[
\Omega \cap \; ]-1,1[^N = \{(\bar x , x_N)\in ]-1,1[^N :\  x_N< 1 -|\bar x|^{\alpha } \}
\]
and that $\partial\Omega$ is $C^2$ outside $]-1,1[^N$; here $\bar x=(x_1, \dots , x_{N-1})$.

Our aim is to obtain stability estimates for the deviation of the eigenvalues and eigenfunctions of the Dirichlet Laplacian $L$ on $\Omega$
from the eigenvalues and eigenfunctions of the Dirichlet Laplacian $L_{\epsilon }$ on the domain $\Omega_{\epsilon }$ defined for $\epsilon\in ]0,1/2[$ by
\begin{eqnarray*}
& &  \Omega_{\epsilon} \cap \; ]-1,1[^N = \{(\bar x , x_N)\in ]-1,1[^N:\ x_N< \min\{1-\epsilon,1 -|\bar x|^{\alpha } \}\}\, \\
& &  \Omega_{\epsilon}\setminus \; ]-1,1[^N = \Omega\setminus ]-1,1[^N \, .
 \end{eqnarray*}

 First we apply Theorem \ref{bootstrap_thm} in the case of the Dirichlet Laplacian on the domain $\Omega $. Let
\[
N_{\alpha} =N+(N-1)\biggl(\frac{1}{\alpha}-1\biggr),
\]
and for a multi-index  $\beta=(\beta_1,\ldots,\beta_N)$, let
\[
|\beta|_{\alpha} =\beta_1 +\cdots +\beta_{N-1} +\alpha\beta_N \, .
\]

\begin{theorem}\label{pcusp}
The Dirichlet Laplacian on $\Omega $ satisfies property (P) for $q_0=\infty$ and any $\gamma>N_{\alpha}/4$.
\end{theorem}
{\em Proof.} We shall apply Theorem \ref{bootstrap_thm}. By \cite[p. 239]{bin} the Sobolev inequality (A) is satisfied with $M=N_{\alpha}$ and
\begin{equation}
\tau(m,\beta,M,p,q) =\frac{|\beta|_{\alpha} + \alpha N_{\alpha}(\frac{1}{p}-\frac{1}{q})}{|\beta|_{\alpha}+\alpha ( m - |\beta|)}.
\label{tau}
\end{equation}
Moreover, by \cite[Theorem 9.1]{mp}, the Dirichlet Laplacian satisfies the a priori estimate (B) for the parameters $m=2$ and arbitrary $1<p_0<\infty $. Hence the result follows by applying (\ref{reg3})  (see Remark~\ref{remboot}). $\hfill\Box$\\

Now let $\epsilon_0\in ]0,1/2[$ be fixed and let $\epsilon\in [0,\epsilon_0 ]$.
We define
\[
\hat \Omega _{\epsilon} =(\Omega\setminus ]-1,1[^N )\cup \{ (\bar x , x_N)\in ]-1,1[^N :\ x_N< h_{\epsilon}(\bar x)\}
\]
where $h_{\epsilon}:]-1,1[^{N-1} \to ]-1,1-\epsilon_0[$ is the locally Lipschitz continuous function implicitly defined by
\begin{equation}
\label{h}
h_{\epsilon}(\bar x ) =1-2\epsilon_0  +\big[ (1-\epsilon_0 -h_{\epsilon}(\bar x))^4   + \max\{|\bar x|^2,\epsilon^{\frac{2}{\alpha }}\}   \big]^{\frac{\alpha }{2}}
\end{equation}
for all $|\bar x |<\epsilon_0^{1/\alpha } $ and by $ h_{\epsilon }(\bar x)=1-|\bar x |^{\alpha }$ for all $|\bar x |\geq
\epsilon_0^{1/\alpha }  $.
Let $\phi_{\epsilon} $ be the map from $\Omega_{\epsilon_0}$ to $\Omega_{\epsilon} $ defined by
\begin{equation}
\label{graf3bis}
\phi_{\epsilon} (\bar x,x_N)\equiv \left\{
\begin{array}{l}
(\bar x,x_N),  \qquad \mbox{ if}\ (\bar x,x_N)\in  \hat \Omega_{\epsilon} ,\\
\Big(\bar x,     -1+2\epsilon_0  +2x_N-\big[ (1-\epsilon_0 -h_{\epsilon }(\bar x))^2 (1-\epsilon_0 -x_N)^2  + \\
\qquad\qquad +\max\{|\bar x|^2,\epsilon^{2/\alpha}\} \big]^{\alpha/2}\Big), \qquad
 \mbox{if}\ (\bar x,x_N)\in \Omega_{\epsilon_0} \setminus\hat \Omega_{\epsilon}.
\end{array}
\right.
\end{equation}
We note that $\phi_{\epsilon} (\Omega_{\epsilon_0} )= \Omega_{\epsilon} $ and $\phi_{\epsilon} (\bar x , h(\bar x ))= (\bar x , h(\bar x )) $, hence $\phi_{\epsilon}\in \Phi (\Omega_{\epsilon_0 } ) $.
Moreover, we note that $\hat \Omega_{\epsilon_0}=\Omega_{\epsilon_0} $ and $\phi_{\epsilon_0}=Id$.

\begin{lemma}
\label{lem:g/g}
Assume that $\frac{1}{2}<\alpha\leq 1$ and $0<\epsilon_0\leq\frac{1}{4}$. There exists a constant $c>0$ depending only on $\alpha $ such that
\begin{equation}
\label{par}
\frac{\det\nabla\phi_{\epsilon}}{\det\nabla\phi_{\epsilon'}}\leq c \; ,
\end{equation}
for all $0\le \epsilon'<\epsilon\le \epsilon_0$.
\end{lemma}
{\em Proof.} We first prove that
\begin{equation}
\label{cl1}
C_{\alpha}\big(\epsilon_0 - \max\{ |\bar x|^{\alpha} ,\epsilon\} \big)
\leq 1-\epsilon_0-h_{\epsilon}(\bar x) \leq \epsilon_0 -  \max\{ |\bar x|^{\alpha} ,\epsilon\},
\end{equation}
for all $0\le \epsilon \le \epsilon_0$, where $ C_{\alpha }=1-1/2^{2\alpha -1}$. Indeed, from (\ref{h}) we have
\begin{eqnarray*}
1-\epsilon_0-h_{\epsilon}&=&\epsilon_0 - \Big[ (1-\epsilon_0-h_{\epsilon})^4 + \max\{ |\bar x|^{2} , \epsilon^{\frac{2}{\alpha}}\} \Big]^{\frac{\alpha}{2}} \\
&\leq& \epsilon_0 - \max\{ |\bar x|^{\alpha} ,\epsilon\},
\end{eqnarray*}
which is the second inequality in (\ref{cl1}). It then follows that
\[
1-\epsilon_0 -h_{\epsilon} \geq \epsilon_0 -\Big[
\big(\epsilon_0 -  \max\{ |\bar x|^{\alpha} ,\epsilon\}\big)^4 + \max\{ |\bar x|^2 ,\epsilon^{\frac{2}{\alpha}}\}   \Big]^{\frac{\alpha}{2}},
\]
hence (\ref{cl1}) will be proved if we show that
\begin{equation}
\label{y}
\epsilon_0 - [ (\epsilon_0 -y)^4 +y^{\frac{2}{\alpha}}]^{\frac{\alpha}{2}} \geq C_{\alpha }(\epsilon_0 -y),
\end{equation}
for all $0<y<\epsilon_0$. To prove (\ref{y}) it suffices to note that since $\epsilon_0\le 1/4$
\[
\epsilon_0 - [ (\epsilon_0 -y)^4 +y^{\frac{2}{\alpha}}]^{\alpha/2} \geq  \epsilon_0  -y -(\epsilon_0 -y)^{2\alpha}
\geq C_{\alpha } (\epsilon_0-y ).
\]
Hence (\ref{cl1}) is proved.

We now prove (\ref{par}). We restrict our attention to $\bar x\in ]-1,1[^{N-1}$ with $|\bar x|\leq\epsilon_0^{1/\alpha}$, since
$\phi_{\epsilon}=\phi_{\epsilon'}$ when
$|\bar x| > \epsilon_0^{1/\alpha}$. We set for simplicity $J=1-\epsilon_0-x_N$.
A direct computation together with (\ref{cl1})  gives
\begin{eqnarray*}
&&\frac{\det\nabla\phi_{\epsilon}}{\det\nabla\phi_{\epsilon'}}=\\
&=&\frac{ 2+ \alpha J\Big[ (1-\epsilon_0- h_{\epsilon})^2 J^2 +
\max\{|\bar x|^2,\epsilon^{\frac{2}{\alpha }}\}\Big]^{ \frac{\alpha -2}{2} } \!\!(1-\epsilon_0-h_{\epsilon})^2 }
{ 2+ \alpha J\Big[ (1-\epsilon_0- h_{{\epsilon '}})^2J^2  +
\max\{|\bar x|^2,{\epsilon '}^{\frac{2}{\alpha }}\}\Big]^{\frac{\alpha-2}{2}} \!\!(1-\epsilon_0-h_{{\epsilon '}})^2 }\\
&\le &C_{\alpha }^{-2}\frac{ 2+ \alpha J\Big[ (\epsilon_0 - \max\{|\bar x|^2,\epsilon^{\frac{2}{\alpha }}\})^2   J^2 +
\max\{|\bar x|^2,\epsilon^{\frac{2}{\alpha }}\}\Big]^{\frac{\alpha-2}{2}} \!\!      (\epsilon_0 - \max\{|\bar x|^2,\epsilon^{\frac{2}{\alpha }}\})^2}
{2+ \alpha J\Big[ (\epsilon_0 - \max\{|\bar x|^2,\epsilon'^{\frac{2}{\alpha }}\})^2   J^2 +
\max\{|\bar x|^2,\epsilon'^{\frac{2}{\alpha }}\}\Big]^{\frac{\alpha-2}{2}} \!\!      (\epsilon_0 - \max\{|\bar x|^2,\epsilon'^{\frac{2}{\alpha }}\})^2} \\
&\le &  C_{\alpha }^{-2}\left\{1+
\frac{
\Big[ (\epsilon_0 - \max\{|\bar x|^2,\epsilon^{\frac{2}{\alpha }}\})^2   J^2 +
\max\{|\bar x|^2,\epsilon'^{\frac{2}{\alpha }}\}\Big]^{\frac{\alpha-2}{2}} \!\!      (\epsilon_0 - \max\{|\bar x|^2,\epsilon ^{\frac{2}{\alpha }}\})^2}
{\Big[ (\epsilon_0-  \max\{|\bar x|^2,\epsilon'^{\frac{2}{\alpha }}\})^2   J^2 +
\max\{|\bar x|^2,\epsilon'^{\frac{2}{\alpha }}\}\Big]^{\frac{\alpha-2}{2}} \!\!      (\epsilon_0 - \max\{|\bar x|^2,\epsilon'^{\frac{2}{\alpha }}\})^2}\right\}\\
&\le  & 2  C_{\alpha }^{-2},
\end{eqnarray*}
where  we have used the simple fact that the function $f(t)=(at+b)^{\frac{\alpha -2}{2}}t$, $t\geq 0$ with $a, b\geq 0$,   is monotone increasing.$\hfill\Box$\\

Recall that $L$, $L_{\epsilon }$ denote the Dirichlet Laplacians on $\Omega$, $\Omega_{\epsilon}$ respectively, as defined in the beginning of this section.

\begin{theorem}
Let $\alpha\in ]1-N/15,1[$. Then the following statements hold.

$\ia$ For all large enough $k\in \N$  and any $\eta >0$ there exists $c>0$ depending only on $\Omega$, $k$ and  $\eta$ such that if $|\Omega\setminus\Omega_{\epsilon }|<c^{-1}$ then the eigenvalues $\lambda_n[L_{\epsilon}]$, $\lambda_n[L]$ satisfy the estimate
\begin{equation}
\bigg(\sum_{n=1}^{\infty}\left| ({\lambda}_n[L_{\epsilon}] +1)^{-k} - ({\lambda}_n[L] +1)^{-k}\right|^2\bigg)^{1/2}
\leq c|\Omega\setminus\Omega_{\epsilon}|^{\frac{1}{2}- \frac{5(1-\alpha )}{N-1+\alpha}   -\eta} .\label{eig}
\end{equation}

$\ib$ Let $\lambda $ be an eigenvalue of multiplicity $m$ of $L$ and let $n\in \N$ be such that
$\lambda  =\lambda_n[L]=\dots =\lambda_{n+m-1}[L]$.
For any $ \eta >0$ there exists $c>0$ depending only on $\Omega $, $\lambda , \lambda_{n-1}[L],  \lambda_{n+m}[L]$ and $\eta$  such that the following is true: if
$|\Omega \setminus\Omega_{\epsilon} |\le c^{-1}$,
then, given orthonormal eigenfunctions
$\psi_n[L_{\epsilon}], \dots , $ $ \psi_{n+m-1}[L_{\epsilon}]$ of $L_{\epsilon}$, there exist corresponding  orthonormal eigenfunctions $ \psi_n[L],$ $ \dots ,$   $  \psi_{n+m-1}[L] $ of $L$ such that
\begin{equation}\label{eigeig}
\| \psi_n[L_{\epsilon}]-\psi_n[L] \|_{L^2(\Omega )}
\le c|\Omega\setminus\Omega_{\epsilon} |^{ \frac{1}{2}- \frac{5(1-\alpha )}{N-1+\alpha}  -\eta}.
\end{equation}
\label{cusp_stability}
\end{theorem}
{\em Proof.} {\it Step 1.} Let $0< \epsilon '<\epsilon \le \epsilon_0 <1/4$. We  apply Theorem~\ref{mainthm} with the maps $\phi = \phi_{\epsilon}:\Omega_{\epsilon _0}\to \Omega_{\epsilon}$ and $\tilde \phi = \phi _{\epsilon' }:\Omega_{\epsilon _0}\to \Omega_{\epsilon' }$. The pull-back to $\Omega_{\epsilon_0}$ of $L_{\epsilon }$ via $\phi_{\epsilon }$  is denoted by $H_{\epsilon }$; similarly, the corresponding matrix $S$ and the function $w$ defined in Section~\ref{perturb} are denoted by $S_{\epsilon ,\epsilon' }$ and $w_{\epsilon , \epsilon'  }$ respectively, and  the operator
$(w_{\epsilon , \epsilon'  }^{-2}\circ \tilde \phi^{(-1)})L_{\epsilon'}$ defined in (\ref{lhat}) is denoted by $\hat L_{\epsilon ,\epsilon' }$;  the matrix $(\nabla\phi_{\epsilon })^{-1}(\nabla\phi_{\epsilon })^{-t}$ is denoted by $a_{\epsilon }$ and the operator $a_{\epsilon }^{1/2}\nabla $ is denoted by $T_{\epsilon}$.  This notation will be used later in Step 3 also for  the limiting case $\epsilon'=0$.

Note that ${\rm det }\nabla \phi_{\epsilon } \geq 1$ and for each $q\in [1, \frac{N}{1-\alpha }[ $ there exists $M>0$ independent of $\epsilon $ such that
\begin{equation}
\label{agiug1}
\|\nabla \phi_{\epsilon}  \|_{L^{q}(\Omega _{\epsilon_0 })}, \| {\rm Adj} (\nabla \phi _{\epsilon}) \|_{L^{q}(\Omega_{\epsilon_0} )}, \|{\rm det }\nabla \phi _{\epsilon } \|_{L^{q}(\Omega_{\epsilon_0 } )}\le M ,
\end{equation}
where ${\rm Adj} (\nabla \phi _{\epsilon})$ denotes the adjugate matrix of $\nabla \phi_{\epsilon} $. Similar computations show that if
$q+1 < \frac{N}{1-\alpha }$ then
\begin{equation}
\label{agiug11}
\|\nabla \phi_{\epsilon}  \|_{L^{q}(\Omega _{\epsilon_0 },g_{\epsilon})},
\|{\rm det }\nabla \phi _{\epsilon } \|_{L^{q}(\Omega_{\epsilon_0 } ,g_{\epsilon})}\le M .
\end{equation}
We now verify that the assumptions of Theorem \ref{mainthm} are satisfied.
It is well-known that $L_{\epsilon }$ and $ L_{\epsilon '}$ satisfy inequality (\ref{est}) with $\alpha =N/2$ and $C_1$ independent of $\epsilon_0$, $\epsilon$ and $\epsilon'$.
Moreover, it follows from \cite[Theorem 3.1]{BS} that there exists $C_1$ independent of $\epsilon_0$, $\epsilon$ and $\epsilon' $ such that (cf. (\ref{lhat}))
$$
\lambda _n[T^*_{\epsilon }S_{\epsilon ,\epsilon'  }T_{\epsilon }] =\lambda _n[\hat L_{\epsilon ,\epsilon '}  ]  \geq C_1n^{\frac{2}{N}},
$$
i.e. $T^*_{\epsilon}S_{\epsilon , \epsilon '}T_{\epsilon }$ satisfies (\ref{est}) with the same parameters.

Now, it is standard that $L_{\epsilon}$ and $ L_{\epsilon' }$ satisfy property (P1) with $q_0=\infty$ and $\gamma =N/4$ (see {\it e.g.}, \cite{babula}).
Since $\psi_n[H_{\epsilon}]=\psi_n[L_{\epsilon}]\circ\phi_{\epsilon}$, $H_{\epsilon}$ also satisfies property (P1) with $q_0=\infty$ and $\gamma =N/4$.
Moreover,  using also (\ref{agiug1}), we have for $q_0$ with $(q_0+2)/2 <N/(1-\alpha)$,
\begin{eqnarray*}
\| \psi_n[ w_{\epsilon,\epsilon'}^{-1}H_{\epsilon'}w_{\epsilon,\epsilon'} ]\|_{L^{q_0}(\Omega_{\epsilon_0},g_{\epsilon})} &=& \| w_{\epsilon,\epsilon'}^{-1}\psi_n[H_{\epsilon'}]\|_{L^{q_0}(\Omega_{\epsilon_0},g_{\epsilon})}\\
&\leq& c\lambda_n[H_{\epsilon'}]^{\frac{N}{4}} \|g_{\epsilon'}^{1/2}\|_{L^{q_0}(\Omega_{\epsilon_0},g_{\epsilon})} \\
&\leq& c\lambda_n[w_{\epsilon,\epsilon'}^{-1}H_{\epsilon'}w_{\epsilon,\epsilon'} ]^{\frac{N}{4}}.
\end{eqnarray*}
Hence the operator $w_{\epsilon,\epsilon'}^{-1}H_{\epsilon'}w_{\epsilon,\epsilon'}$ satisfies property (P1) for any $q_0<2(N-1+\alpha)/(1-\alpha)$ and $\gamma=N/4$, uniformly in $\epsilon_0,\epsilon,\epsilon'$.

By the argument in \cite[Theorem 9.1]{mp} (which deals with the case of a cusp) it follows that
the operator $L_{\epsilon }$  satisfies the {\em a priori} estimate (B) with any $p_0>1$, $m=2$ and  $A_{p}$ independent of $\epsilon $. Since the Sobolev inequality (A) is also valid  with $M=N_{\alpha }$ and $\tau $ defined by (\ref{tau}), and  $C_4$ independent of $\epsilon$, by Theorem~\ref{bootstrap_thm} it follows that the operator $L_{\epsilon }$  satisfies property (P2) for $q_0=\infty $ and any $\gamma >N_{\alpha }/4$, uniformly in $\epsilon$ (see also Theorem~\ref{pcusp}).
Since $\nabla\psi_n[H_{\epsilon}]= (\nabla\psi_n[L_{\epsilon}]\circ\phi_{\epsilon})\nabla\phi_{\epsilon}$, we have for any $q_0$ with
$q_0+1<\frac{N}{1-\alpha }$ (cf. (\ref{agiug11})) and any $\eta>0$,
\begin{eqnarray*}
\|\nabla\psi_n[H_{\epsilon}]\|_{L^{q_0}(\Omega_{\epsilon_0},g_{\epsilon})} &\leq& c\lambda_n[H_{\epsilon}]^{\half +\frac{N_{\alpha}}{4}+\eta}
\|\nabla \phi_{\epsilon}  \|_{L^{q_0}(\Omega _{\epsilon_0 },g_{\epsilon})}\\
&\leq & c\lambda_n[H_{\epsilon}]^{\half +\frac{N_{\alpha}}{4}+\eta},
\end{eqnarray*}
uniformly in $\epsilon,\epsilon_0$. Hence $H_{\epsilon}$ satisfies property (P2) for any $q_0<(N-1+\alpha)/(1-\alpha)$ and any $\gamma>N_{\alpha}/4$.

We finally consider $T_{\epsilon}^*S_{\epsilon,\epsilon'}T_{\epsilon}$. By Lemma~\ref{lem:g/g}, $w_{\epsilon ,\epsilon'}^2\le c$ hence the operator
$\hat L_{\epsilon ,\epsilon '}=(w^{-2}_{\epsilon ,\epsilon'}  \circ\phi_{\epsilon'}^{(-1)})L_{\epsilon'}$, which is self-adjoint on
$L^{2} (\Omega_{\epsilon '},w_{\epsilon,\epsilon'}^2\circ \phi_{\epsilon'}^{(-1)} )$,
also satisfies the {\em a priori} estimate (B), for the same parameters as
$L_{\epsilon'}$. Since the Sobolev inequality (A) is also valid (cf. Theorem~\ref{pcusp}), we can apply Theorem \ref{bootstrap_thm} and (\ref{reg3}) and conclude that any  eigenfunction $\psi_n[\hat L_{\epsilon , \epsilon '}]$ of $\hat L_{\epsilon , \epsilon '}$  satisfies
\begin{equation}
\label{apunw}
\| D^{\beta }\psi_n[\hat L_{\epsilon , \epsilon '}]\|_{L^{\infty}(\Omega_{\epsilon '})} \leq
c \lambda_n[\hat L_{\epsilon , \epsilon '}]^{\frac{|\beta  |}{2}+\frac{N_{\alpha}}{2p_0}+\eta}\| \psi_n[\hat L_{\epsilon , \epsilon '}] \|_{L^{p_0} (\Omega_{\epsilon '} )} \, ,
\end{equation}
for all multi-indeces $\beta $ with $|\beta |\le 1$, all $p_0>1$ and any $\eta>0$, uniformly in $\epsilon_0,\epsilon,\epsilon'$.
Now, for any    $p_0$ with $1<p_0< 2(1-(1-\alpha)/N)$ we have by H\"{o}lder inequality,
\begin{eqnarray}
\| \psi_n[\hat L_{\epsilon , \epsilon '}] \|_{L^{p_0} (\Omega_{\epsilon '} )}&\leq&
\|\psi_n[\hat L_{\epsilon , \epsilon '}](w_{\epsilon,\epsilon'}\circ\phi_{\epsilon'}^{(-1)}) \|_{L^{2} (\Omega_{\epsilon '} )}\|w_{\epsilon,\epsilon'}^{-1}\circ\phi_{\epsilon'}^{(-1)}\|_{L^{\frac{2p_0}{2-p_0}} (\Omega_{\epsilon '} )}\nonumber\\
&=&\|w_{\epsilon,\epsilon'}^{-1}\circ\phi_{\epsilon'}^{(-1)}\|_{L^{\frac{2p_0}{2-p_0}} (\Omega_{\epsilon '} )}\nonumber\\
&\leq&\|g_{\epsilon'}^{1/2}\circ\phi_{\epsilon'}^{(-1)}\|_{L^{\frac{2p_0}{2-p_0}} (\Omega_{\epsilon '} )}\label{p_0}\\
&=&\bigg(\int_{\Omega_{\epsilon_0}}g_{\epsilon'}^{\frac{2}{2-p_0}}dx\bigg)^{\frac{2-p_0}{2p_0}}\nonumber\\
&\leq& c , \nonumber
\end{eqnarray}
uniformly in $\epsilon_0,\epsilon,\epsilon'$. Now, we have $\psi_n[T^*_{\epsilon}S_{\epsilon , \epsilon '}T_{\epsilon }]=
\psi_n[\hat L_{\epsilon , \epsilon '}]\circ\phi_{\epsilon'}$. Hence (\ref{apunw}) implies that $T^*_{\epsilon}S_{\epsilon , \epsilon '}T_{\epsilon }$ satisfies (P1) for $q_0=\infty$ and any $\gamma>N_{\alpha}N/(4(N-1+\alpha))$. Moreover, for any $p_0$ as in (\ref{p_0}), any $\eta>0$ and any $q_0$ with $q_0+1<N/(1-\alpha)$ we have, using also (\ref{agiug11}),
\begin{eqnarray*}
\|\nabla\psi_n[T^*_{\epsilon}S_{\epsilon , \epsilon '}T_{\epsilon }]\|_{L^{q_0}(\Omega_{\epsilon_0},g_{\epsilon})}&\leq&
c \lambda_n[T^*_{\epsilon}S_{\epsilon , \epsilon '}T_{\epsilon }]^{\frac{1}{2}+\frac{N_{\alpha}}{2p_0}+\eta}
\|\nabla\phi_{\epsilon'}\|_{L^{q_0}(\Omega_{\epsilon_0},g_{\epsilon})}\\
&\leq& c\lambda_n[T^*_{\epsilon}S_{\epsilon , \epsilon '}T_{\epsilon }]^{\frac{1}{2}+\frac{N_{\alpha}}{2p_0}+\eta}.
\end{eqnarray*}
Hence $T^*_{\epsilon}S_{\epsilon , \epsilon '}T_{\epsilon }$ satisfies property (P2) for any $q_0<(N-1+\alpha)/(1-\alpha)$ and any $\gamma>N_{\alpha}N/(4(N-1+\alpha))$.

Summing up, Theorem \ref{mainthm} can be applied for any $q_0<(N-1+\alpha)/(1-\alpha)$ and any $\gamma>N_{\alpha}N/(4(N-1+\alpha))$.

Applying the theorem we obtain that for any $q_0<(N-1+\alpha)/(1-\alpha)$ and $k\in {\N}$ sufficiently large there holds
\begin{equation}
\|(w^{-1}_{\epsilon ,\epsilon '}H_{\epsilon '}w_{\epsilon , \epsilon'} +1)^{-k} -(H_{\epsilon }+1 )^{-k}\|_{\cC^2(L^2(\Omega_{\epsilon_0}),g_{\epsilon})}\leq c\delta _{\frac{2q_0}{q_0-2}}(\phi_{\epsilon }, \phi_{\epsilon '}),
\label{res22}
\end{equation}
uniformly in $\epsilon_0 $, $\epsilon $, $\epsilon'$, provided $\delta _{\frac{2q_0}{q_0-2}}(\phi_{\epsilon }, \phi_{\epsilon '})<c^{-1}$.
Since $w_{\epsilon_0,\epsilon'}=w_{\epsilon,\epsilon'}w_{\epsilon_0,\epsilon}$, by unitary equivalence and (\ref{res22}) we obtain
\begin{equation}
\|(w^{-1}_{\epsilon_0 ,\epsilon '}H_{\epsilon '}w_{\epsilon_0 , \epsilon'} +1)^{-k} -
(w^{-1}_{\epsilon_0 ,\epsilon}H_{\epsilon}w_{\epsilon_0 , \epsilon} +1)^{-k}\|_{\cC^2(L^2(\Omega_{\epsilon_0}))}\leq c\delta _{\frac{2q_0}{q_0-2}}(\phi_{\epsilon }, \phi_{\epsilon '}),
\label{res22a}
\end{equation}
uniformly in $\epsilon_0 $, $\epsilon $, $\epsilon'$. In particular, for $\epsilon =\epsilon_0$
\begin{equation}
\|(w^{-1}_{\epsilon_0 ,\epsilon '}H_{\epsilon '}w_{\epsilon_0 , \epsilon'} +1)^{-k} -
(H_{\epsilon_0} +1)^{-k}\|_{\cC^2(L^2(\Omega_{\epsilon_0}))}\leq c\delta _{\frac{2q_0}{q_0-2}}(\phi_{\epsilon_0 }, \phi_{\epsilon '}).
\label{res22a''}
\end{equation}

{\em Step 2.} We now estimate the right-hand side of (\ref{res22a}). We note that
\begin{eqnarray}\label{agiug2}
& & |S^{1/2}_{\epsilon , \epsilon '}|\le  c |\nabla \phi_{\epsilon }|\, |{\rm Adj}(\nabla \phi_{\epsilon' } )  |,\nonumber \\
& &  |S^{-1/2}_{\epsilon ,\epsilon '}|\le  c  |\nabla \phi_{\epsilon '}|\, |{\rm Adj}(\nabla \phi_{\epsilon } )  |,\nonumber \\
& & |a_{\epsilon }^{1/2}|\le c |{\rm Adj}(\nabla \phi_{\epsilon } ) |
\end{eqnarray}
for some constant $c>0$. Note also that $\phi_{\epsilon }=\phi_{\epsilon '}$ on $U_{\epsilon}=\{(\bar x,x_N)\in \Omega_{\epsilon_0}:\ |\bar x|>\epsilon ^{1/\alpha } \}$.
Thus, by  (\ref{agiug1}), (\ref{agiug2}) and  H\"{o}lder inequality it follows that  if $1<s<q_0/6<q_0<\frac{N-1+\alpha}{1-\alpha}$
\begin{eqnarray*}\lefteqn{
\| (S^{1/2}_{\epsilon ,\epsilon' }-S^{-1/2}_{\epsilon ,\epsilon ' })a_{\epsilon }^{1/2}
 \|_{L^s(\Omega_{\epsilon_0 } ,g_{\epsilon})} }\nonumber \\
& &
\le  | \Omega _{\epsilon_0} \setminus U_{\epsilon} |^{\frac{1}{s}-\frac{4}{q_0}} \| (S^{1/2}_{\epsilon ,\epsilon' }-S^{-1/2}_{\epsilon ,\epsilon ' })a_{\epsilon }^{1/2}g_{\epsilon}^{1/s} \|_{L^{q_0/4}(\Omega _{\epsilon_0} \setminus U_{\epsilon}  )}
\le c | \Omega _{\epsilon_0} \setminus U_{\epsilon} |^{\frac{1}{s}-\frac{4}{q_0}},
\end{eqnarray*}
and
\begin{eqnarray*}\lefteqn{
\| (S_{\epsilon ,\epsilon' }-I)a_{\epsilon }^{1/2} \|_{L^s(\Omega_{\epsilon_0 },g_{\epsilon} )} }\nonumber \\
& &
\le  |\Omega _{\epsilon_0} \setminus U_{\epsilon} |^{\frac{1}{s}-\frac{6}{q_0}} \| (S_{\epsilon ,\epsilon' }-I)a_{\epsilon }^{1/2} g_{\epsilon}^{1/s} \|_{L^{q_0/6}(\Omega _{\epsilon_0} \setminus U_{\epsilon} )}
\le c | \Omega _{\epsilon_0} \setminus U_{\epsilon}  |^{\frac{1}{s}-\frac{6}{q_0}},
\end{eqnarray*}
for some constant $c>0$.
One can similarly estimate the other summands in (\ref{vicinity}) and  get
\begin{equation}
\delta_s(\phi_{\epsilon },\phi_{\epsilon '} )\le c |\Omega_{\epsilon_0}\setminus U_{\epsilon }|^{\frac{1}{s}-\frac{6}{q_0}} \; .
\label{nat}
\end{equation}
uniformly in $\epsilon_0,\epsilon$ and $\epsilon'$, for $s$ and $q_0$ as above.

{\it Step 3.}
Since $\alpha > 1-N/15$, it is possible to choose $14<q_0<(N-1+\alpha)/(1-\alpha )$ which guarantees that $2q_0/(q_0-2)< q_0/6 $; thus choosing $s=2q_0/(q_0-2)$
in (\ref{nat}) it follows in particular that $\delta_s(\phi_{\epsilon },\phi_{\epsilon '} )\to 0$ as $\epsilon \to 0$, uniformly in $\epsilon'\in (0,\epsilon)$. This combined with (\ref{res22a}) implies that
the sequence $(w^{-1}_{\epsilon_0 ,\epsilon }H_{\epsilon }w_{\epsilon_0 , \epsilon} +1)^{-k}$ is Cauchy in
$\epsilon$ for $\epsilon \to 0$.
Thus, by passing to the limit in (\ref{res22a''})
as $\epsilon '\to 0$ we obtain
\begin{equation}
\|(w^{-1}_{\epsilon_0 ,0 }H_{0}w_{\epsilon_0 , 0} +1)^{-k} -(H_{\epsilon_0}+1 )^{-k}\|_{\cC^2(L^2(\Omega_{\epsilon_0}))}
\leq c\delta _{\frac{2q_0}{q_0-2}}(\phi_{\epsilon_0 }, \phi_0).
\label{res222}
\end{equation}
Taking into account that $\phi_{\epsilon_0}=Id$, using (\ref{agiug2}) and proceeding as in Step 2, we get
\begin{equation}
\label{nat2}
\delta _{s}(\phi_{\epsilon_0 }, \phi_0)\le c | \Omega_{\epsilon_0}\setminus \hat\Omega_{0 }   |^{\frac{1}{s}-\frac{4}{q_0}}
\end{equation}
for all $1<s<q_0/4<q_0<(N-1+\alpha)/(1-\alpha )$.
Since $1-\epsilon_0 -h_{0} \le \epsilon_0 -|\bar x |^{\alpha }$ ({\it cf.} (\ref{cl1}) with $\epsilon =0$), we get  $|\Omega_{\epsilon_0 } \setminus \hat\Omega_{0}  |\le | \Omega  \setminus \Omega_{\epsilon_0} |$.
By (\ref{res222}), (\ref{nat2}) and choosing $s=2q_0/(q_0-2)$ it follows that
\begin{equation}
\|(w^{-1}_{\epsilon_0 ,0 }H_{0}w_{\epsilon_0 , 0} +1)^{-k} -
(H_{\epsilon_0}+1 )^{-k}\|_{\cC^2(L^2(\Omega_{\epsilon_0}))}\leq c  | \Omega  \setminus \Omega_{\epsilon_0} |^{\frac{q_0-10}{2q_0}}.
\end{equation}
In order to conclude, it suffices to observe that $(q_0-10)/(2q_0)\to \frac{1}{2}- \frac{5(1-\alpha )}{N-1+\alpha}$
as $q_0\to (N-1+\alpha)/(1-\alpha)$ and proceed as in the proof of Theorems~\ref{thm:series} and \ref{thm:eigen}. \hfill $\Box$

{\small
\noindent Gerassimos Barbatis\\
Department of Mathematics\\
University of Athens\\
157 84 Athens\\
Greece\\
e-mail: gbarbatis@math.uoa.gr\\

\noindent Pier Domenico lamberti\\
Dipartimento di matematica Pura ed Applicata\\
Universit\`{a} degli Studi di Padova\\
Via Trieste, 63\\
35121 Padova\\
Italy\\
e-mail: lamberti@math.unipd.it
}

\end{document}